\documentclass[12pt]{article}
\usepackage{cite,epsfig}
\usepackage{graphicx}
\usepackage{amsthm,amsmath,amssymb}
\newcommand{\rstate}[1]{\left| #1 \right>}

\newcommand{\ba}{\begin{array}}
\newcommand{\ea}{\end{array}}
\newcommand{\be}{\begin{equation}}
\newcommand{\ee}{\end{equation}}

\newcommand{\rad}{\operatorname{rad}}

\newcommand{\bin}[2]{\left( \begin{array}{@{}c@{}} #1 \\ #2 \end{array} \right)}

\def\half{\frac{1}{2}}

\def\f{\frac}

\def\TL{Temperley-Lieb }

\def\C{\mathbb{C}}
\def\Z{\mathbb{Z}}

\newtheorem{thm}{Theorem}[section]
\newtheorem{defn}[thm]{Definition}
\newtheorem{prop}[thm]{Proposition}
\newtheorem{lemma}[thm]{Lemma}
\newtheorem{remark}[thm]{Remark}
\newtheorem{corol}[thm]{Corollary}
\newtheorem{example}[thm]{Example}
\newtheorem{conj}[thm]{Conjecture}

\begin{document}
\begin{titlepage}

\begin{flushright}

\end{flushright}

\vskip 2 cm

\begin{center}
{\LARGE The two-boundary \TL algebra}
\vskip 1 cm

{\large Jan de
  Gier\footnote{degier@ms.unimelb.edu.au} and Alexander  Nichols\footnote{nichols@sissa.it} }

\begin{center}
{\em $^1$Department of Mathematics and Statistics, \\
University of Melbourne, \\ Victoria VIC 3010, Australia.}

{\em $^2$International School for Advanced Studies, \\
Via Beirut 1, 34100 Trieste, Italy.}

{\em$^2$ Istituto Nazionale di Fisica Nucleare (INFN), \\
Sezione di Trieste, Trieste, Italy.}

\vskip 1 cm
\end{center}

\vskip .5 cm

\begin{abstract}
We study a two-boundary extension of the \TL algebra which has
recently arisen in  statistical mechanics. This algebra lies in a
quotient of the affine Hecke algebra of type C and has a natural
diagrammatic representation. The algebra has three parameters and, for
generic values of these, we determine its representation theory. 

We use the action of the centre of the affine Hecke algebra to show
that all irreducible representations lie within a finite
dimensional diagrammatic quotient. These representations are fully
characterised by an additional parameter related to the action of the
centre. For generic values of this parameter there is a unique
representation of dimension $2^N$ and we show that it is isomorphic to a
tensor space representation. We construct a basis in which the Gram
matrix is diagonal and use this to discuss the irreducibility of this
representation. 
\end{abstract}

\end{center}

\end{titlepage}
%
%
%
\section{Introduction}
The \TL (TL) algebra \cite{TL71} first appeared in statistical mechanics as a tool
to analyze various interrelated lattice models such as the
$Q$-state Potts model, the O($n$) loop model and the six-vertex model,
see e.g. \cite{B82,M91}. It subsequently played a
crucial role in both mathematics and in physics, for example in the construction of
knot invariants \cite{J91} and the development of solvable lattice
models \cite{B82}. The O($n$) loop models, which appear in the
diagrammatic representation of the TL algebra \cite{K87}, have
attracted renewed attention recently in the context of stochastic
Loewner evolution (SLE) \cite{S00} - see e.g. \cite{KN04} for a review. 

In this paper we shall study a two-boundary extension of the TL
algebra which naturally arises from considering the addition of
integrable boundary terms to the six-vertex model
\cite{deGP04,deGNPR05,N06a,N06b}. This algebra first appeared in
\cite{MNGB04} where the O($1$) model was studied, see also
\cite{deG05}. This model describes critical bond percolation and is
equivalent to the stochastic raise and peel model
\cite{deGNPR04}. The two-boundary TL algebra also underlies the
partially asymmetric exclusion process with open boundaries
\cite{deGE05}. It generalises the one-boundary Temperley-Lieb (1BTL), or
blob \cite{MS93,MS94,MW00,MW03}, algebra by the addition of a second
boundary generator.  

An important open problem in the theory of solvable lattice models is
the construction of so-called Bethe Ansatz equations for the
six-vertex model with general integrable boundary terms. Recent progress was made in
\cite{CLSW03,N02} -- see \cite{deGP04} for a loop model context -- where this construction
was achieved for certain special cases. In \cite{deGNPR05} we were led
to the conjecture that these cases were related to properties of the 2BTL representation theory. This paper largely arose
from the desire to understand better the 2BTL algebra and its
representation theory.  

In stark contrast to the TL and 1BTL algebras the 2BTL algebra,
defined in Section \ref{sec:TLboundaryextns}, is infinite
dimensional. This fact makes the study of its representation theory
considerably more interesting. The algebra contains three parameters
and for generic values of these we determine its irreducible
representations.  

As is the case for the TL and 1BTL algebras, the 2BTL algebra has a simple diagrammatic
representation which we give in Section \ref{sec:DiagrammaticRepn}. The 2BTL algebra
is a quotient of the affine Hecke algebra of type C. This algebra has
a large centre which is conveniently described using a commutative set
of Murphy elements. We review this in Section \ref{sec:Hecke}. In Theorem
\ref{thm:IJI=bI} we show, using Schur's lemma, that all irreducible
representations of the 2BTL algebra lie within simple diagrammatic
quotients. 

The irreducible representations are fully characterised by an additional
parameter $b$ related to the action of the centre. For generic values of
this parameter we find a unique largest irreducible representation, called
$W^{(N)}(b)$, of dimension $2^N$. We study the structure of this
representation by constructing a basis ${\bf B}_1$ which diagonalizes
the Murphy elements of the type B Hecke algebra. The type B Hecke algebra is
related to the 1BTL algebra, but in the basis ${\bf B}_1$ all the
generators, including \emph{both} boundary generators, act in a simple
way. We prove the Gram matrix is diagonal in basis ${\bf B}_1$ and in
Theorem \ref {thm:GramDet1} we compute its determinant. For generic
values of $b$ the representation $W^{(N)}(b)$ is irreducible however
it fails to be at a discrete set of points. 

In Section \ref{sec:SpinChainRepn} we show that the $2^N$ dimensional
representation, $W^{(N)}(b)$, is isomorphic to a tensor product
representation. In this representation all the parameters acquire a
physical significance. The points where the action of the 2BTL
generators becomes indecomposable are exactly those previously
conjectured in \cite{deGNPR05}. 

In Section \ref{sec:OtherIrreps} we discuss the cases in which the
action of the 2BTL generators in the representation $W^{(N)}(b)$
becomes reducible but indecomposable. At these points the action of
the centre takes only a discrete set of values and it is possible to
construct smaller irreducible representations. A large number of these
are found to have a simple diagrammatic description. 
\section*{Acknowledgements}
A.~N. would like to thank B.~Westbury for collaboration during the 
early stages of this work. It is a pleasure to thank V.~Rittenberg and 
P.~Pyatov for discussions. We would also like to thank A.~Ram for
useful conversations and for emphasizing to us the crucial role of the
Murphy elements in the theory of Hecke, and affine Hecke, algebras. 

\section{Definition of algebras}

We will start by defining the main algebras which we will study in
this paper. These algebras are all quotients of Hecke algebras which
will be described in Section \ref{sec:Hecke}. 

\subsection{Boundary extensions of the Temperley-Lieb algebra}
\label{sec:TLboundaryextns}
\begin{defn}
\label{defn:TL}
Let $\delta$ be an indeterminate. The \TL (TL) algebra \cite{TL71,M91} is an associative algebra, over the ring $\Z[\delta]$, given by generators $e_i$ with $i=1,\cdots,N-1$ obeying the relations: 
\begin{align*}
e_i^2 &= \delta ~e_i, \\
e_i e_{i\pm 1} e_i &=e_i, \\
e_i e_j &= e_j e_i \qquad  |  i-j |  >1. 
\end{align*}
\end{defn}
\begin{defn}
\label{defn:1BTL}
Let $s_1$ and $\delta$ be indeterminates. The one-boundary \TL (1BTL) algebra \cite{MS93,MS94,MW00,MW03} is an associative algebra over $\Z[\delta,s_1]$ defined by adding an additional generator $e_0$ to the TL algebra. This generator is required to satisfy: 
\begin{align*} 
e_0^2&=s_1 ~e_0, \\
e_1 e_0 e_1&=e_1, \\
e_i e_0 & = e_0 e_i \qquad i>1.
\end{align*}
\end{defn}
There is also an analogous (isomorphic) algebra defined by instead
adjoining a generator $e_N$ at the right hand end. 
The primary object of study in this paper is the algebra with boundary
generators at \emph{both} ends. 
\begin{defn}
\label{defn:2BTL}
Let $s_1$, $s_2$ and $\delta$ be indeterminates. The two-boundary Temperley-Lieb (2BTL) algebra is an associative algebra over $\Z[\delta,s_1,s_2]$ defined by adding an additional generator $e_N$ to the 1BTL algebra. This generator is required to satisfy: 
\begin{align*}
e_{N}^2&=s_2 ~e_{N}, \\
e_{N-1} e_{N} e_{N-1}&=e_{N-1}, \\
e_i e_{N} &= e_{N} e_i \qquad i<N-1.
\end{align*}
\end{defn}
The 2BTL algebra appeared in recent physics literature and its representation theory is important in applications such as the quantum XXZ spin chain and O($n=1$) loop model, both with open boundaries, and the combinatorics of alternating sign matrices \cite{MNGB04,deG05,deGP04,deGNPR05}, as well as for the partially asymmetric exclusion process \cite{deGE05}.

Each of the algebras TL, 1BTL and 2BTL arises as a quotient of a Hecke
algebra. This fact will be formulated more precisely in
Proposition~\ref{prop:HeckeQuotients} below after we have first defined the
relevant Hecke algebras.

\subsection{Relevant Hecke algebras}
\label{sec:Hecke}
\subsubsection{Hecke algebra of type A}
\begin{defn} Let $q\in \C^*$ be an indeterminate. The Hecke algebra of type A \cite{J87,DJ86,DJ87} over $\Z[q,q^{-1}]$ has generators $g_i$ with $i=1,2,\cdots,N-1$ obeying the relations: 
\begin{align*}
g_i g_{i+1} g_i &=  g_{i+1} g_i g_{i+1} \qquad \quad i=1,2,\cdots,N-2 \\
g_i g_j &= g_j g_i  \qquad\qquad \qquad |i-j|>1\\
\left( g_i -q \right) \left( g_i + q^{-1} \right) &=0
\end{align*}
We also define the Murphy elements of type A:
\begin{align*}
J_1^{(\rm A)}&= g_1^2 \\
J_{i}^{(\rm A)}&= g_{i} J_{i-1}^{(\rm A)} g_{i} \quad \quad 2 \le i \le N-1
\end{align*}
\end{defn}
\begin{prop}
\label{prop:MurphyTypeA}
The Murphy elements $J_i^{(\rm A)}$ are pairwise
commuting and obey the following relations: 
\begin{align*} 
\left[g_1,J_j^{(\rm A)}\right]&=0 \qquad j\geq 1,\\
\left[g_i,J_j^{(\rm A)}\right]&=0 \qquad i \ge 2,~j \ne i-1,i \\
\left[g_i,J_i^{(\rm A)} J_{i-1}^{(\rm A)}\right]&=0 \qquad i \ge 2, \\
\left[g_i,J_i^{(\rm A)} + J_{i-1}^{(\rm A)}\right]&=0 \qquad i \ge 2. 
\end{align*}
These statements imply that all completely symmetric polynomials in the set
$\{ J_i^{(\rm A)} \}$ are central. 
\end{prop}
\begin{proof}
These statements are all simple to prove and here we shall only present proofs of the final
two. For $i \ge 2$ we have, omitting for brevity the superscript $(\rm
A)$:
\begin{equation*}
[g_i,J_{i-1}J_i]= [g_i,J_{i-1} g_i J_{i-1} g_i ] = [J_i,J_{i-1}]g_i =0.
\end{equation*}
\begin{equation*}
[g_i,J_{i-1}+J_i]=[g_i,J_{i-1} + g_i J_{i-1} g_i ] = g_i J_{i-1}-J_{i-1}g_i + g_i^2 J_{i-1} g_i - g_i J_{i-1} g_i^2 =0
\end{equation*}
where in the second line we have used the Hecke relation $g_i^2 = (q-q^{-1})g_i + 1$. Note that in proving these results we only need commutativity of the
Murphy elements, their inductive definition, and the Hecke condition on the generators.
We shall therefore be able to repeat this argument 
later with modified expressions for the first Murphy element. 
\end{proof}
\subsubsection{Hecke algebra of type B}
\begin{defn}
\label{defn:HeckeB}
Let $q^{\omega_1}\in\C$ be an indeterminate. The Hecke algebra of type B over $\Z[q,q^{-1},q^{\omega_1},q^{-\omega_1}]$ \cite{DJ86,DJ87,DJ92,DJM95} is given by adjoining to the Hecke algebra of type A the additional generator $g_0$. We have the relations: 
\begin{align*}
g_0 g_1 g_0 g_1 &= g_1 g_0 g_1 g_0 \\
g_0 g_i &= g_i g_0 \quad \quad i>1  \\
\left( g_0 - q^{\omega_1} \right) \left( g_0 - q^{-\omega_1} \right) &=0
\end{align*}
We also define the Murphy elements of type B:
\begin{align*}
J_0^{(\rm B)}&= g_0  \\
J_{i}^{(\rm B)}&= g_{i} J_{i-1}^{(\rm B)} g_{i} \quad \quad 1 \le i \le N-1
\end{align*}
\end{defn}
\begin{prop}
\label{prop:MurphyTypeB}
The Murphy elements $J_i^{(\rm B)}$ are pairwise commuting. For $i \ge 1$ we have the following relations:
\begin{align*} 
\left[g_i,J_j^{(\rm B)} \right]&=0 \quad \quad j \ne i-1,i \\
\left[g_i,J_{i-1}^{(\rm B)} J_i^{(\rm B)} \right]&=0  \\
\left[g_i,J_{i-1}^{(\rm B)} + J_i^{(\rm B)} \right]&= 0 
\end{align*}
These statements imply that all completely symmetric polynomials in the set
$\{ J_i^{(\rm B)} \}$ are central.
\end{prop}
\begin{proof}
Similar to the proof of Proposition~\ref{prop:MurphyTypeA}.
\end{proof}
\subsubsection{Affine Hecke algebra of type C}
\begin{defn}
\label{defn:AffineTypeC}
Let $q^{\omega_2}\in\C$ be an indeterminate. The affine Hecke algebra of type C over $\Z[q,q^{-1},q^{\omega_1},q^{-\omega_1},q^{\omega_2},q^{-\omega_2}]$ \cite{L89,R03,RR03} is given by adding to the Hecke algebra of type B an additional boundary generator $g_N$ with relations 
\begin{align*}
g_N g_{N-1} g_N g_{N-1} &= g_{N-1} g_N g_{N-1} g_N \\
g_N g_{i} &= g_{i} g_N \quad \quad \quad \quad 0 \le i\le N-2 \\
\left( g_N - q^{\omega_2} \right) \left( g_N - q^{-\omega_2} \right) &=0
\end{align*}
We also define the Murphy elements for the affine Hecke algebra of
type C:
\begin{align*}
J^{(\rm C)}_0&= g_1^{-1} g_2^{-1} \cdots g_{N-1}^{-1} g_N g_{N-1} \cdots g_2 g_1 g_0\\
J^{(\rm C)}_{i}&= g_{i} J^{(\rm C)}_{i-1} g_{i} \quad \quad 1 \le i \le N-1.
\end{align*}
\end{defn}
The inductive definition of the Murphy elements is the same as the
previous cases. We shall also use another equivalent definition of the
affine Hecke algebra of type C, where instead of $g_N$ we take
$J_0^{(\rm C)}$ to be the additional independent generator: 
\begin{prop}
\label{prop:EquivalentAffineC}
Let $q^{\omega_2}\in\C$ be an indeterminate. The affine Hecke algebra of type C can be equivalently described by
adding to the Hecke algebra of type B an additional generator $J^{(\rm C)}_0$
with relations: 
\begin{align}
\label{eqn:EquivalentAffineC}
g_i J^{(\rm C)}_0 &= J^{(\rm C)}_0 g_i \quad \quad i>1 \nonumber\\
J_0^{(\rm C)} g_1 J^{(\rm C)}_0 g_1 &= g_1 J^{(\rm C)}_0 g_1 J_0^{(\rm C)}\\
g_0 g_1 J^{(\rm C)}_0 g_1 &= g_1 J^{(\rm C)}_0 g_1 g_0 \nonumber \\
\left( J^{(\rm C)}_0 g_0^{-1} - q^{\omega_2} \right) \left( J^{(\rm C)}_0 g_0^{-1} - q^{-\omega_2} \right) &= 0 \nonumber
\end{align}
The set of Murphy elements are defined inductively from $J^{(\rm C)}_0$ as in Definition \ref{defn:AffineTypeC}.
\end{prop}
\begin{proof}
The relations \eqref{eqn:EquivalentAffineC} follow by simple algebra from Definition \ref{defn:AffineTypeC}. Conversely defining $g_N$ by:
\begin{equation}
\label{eqn:gNExpression}
g_N=g_{N-1} g_{N-2} \cdots g_{1} J^{(\rm C)}_0 g_0^{-1} g_1^{-1} \cdots g_{N-1}^{-1}
\end{equation}
and using the relations given in \eqref{eqn:EquivalentAffineC} we recover those of Definition \ref{defn:AffineTypeC}.
\end{proof}
\begin{prop}
\label{prop:AffineTypeC}
The Murphy elements $J^{(\rm C)}_i$ are pairwise commuting and obey the following relations:
\begin{align}
\label{eqn:AffineTypeCcommutation}
\left[g_0,J^{(\rm C)}_j \right]&=0 \qquad j \ne 0, \nonumber \\
\left[g_i,J^{(\rm C)}_j \right]&=0 \qquad 1 \le i \le N-1, j \ne i-1,i, \nonumber\\
\left[g_i,J^{(\rm C)}_{i-1} J^{(\rm C)}_i \right]&=0 \qquad 1 \le i \le N-1, \\
\left[g_i,J^{(\rm C)}_{i-1} + J^{(\rm C)}_i \right]&=0 \qquad 1 \le i \le N-1, \nonumber \\
\left[g_0,J^{(\rm C)}_0+\left(J^{(\rm C)}_0\right)^{-1} \right]&=0. \nonumber
\end{align}
These statements imply that all completely symmetric
polynomials in $\{J^{(\rm C)}_i,\left(J^{(\rm C)}_i\right)^{-1}\}$ are central. 
\end{prop}
\begin{proof}
The proofs of all but the final equation in
\eqref{eqn:AffineTypeCcommutation} are the same as in 
Propositions~\ref{prop:MurphyTypeA} and \ref{prop:MurphyTypeB}. We
shall therefore not repeat them here. To prove the last equation,
omitting the superscript (C), we note that: 
\[
J_0^{-1} = \left(q^{\omega_2} + q^{-\omega_2} \right) g_0^{-1} - g_0^{-1} J_0 g_0^{-1},
\]
so that
\[ 
[g_0,J_0+J_0^{-1}] = g_0 J_0 - J_0 g_0 - J_0 g_0^{-1} + g_0^{-1}
J_0 = [g_0+g_0^{-1},J_0] = 0,
\]
as by Definition \ref{defn:HeckeB} we have $g_0+g_0^{-1}= q^{\omega_1} + q^{-\omega_1} $.

Now using equations \eqref{eqn:AffineTypeCcommutation} we find that
all completely symmetric polynomials in the the set $\{J_i,J_i^{-1}
\}$ commute with $g_0$ and $g_i$ for $1 \le i \le N-1$. The
commutation with $g_N$ follows using \eqref{eqn:gNExpression}.
\end{proof}
\begin{remark}
\label{remark:NogNgenerator}
It is only the generators given in Proposition \ref{prop:EquivalentAffineC} which appear in \eqref{eqn:AffineTypeCcommutation}. We did not find any simple relations between the Murphy elements of the affine type C Hecke algebra, $J^{(\rm C)}_{i}$, and the generator $g_N$.
\end{remark}
\subsubsection{\TL quotients}

The \TL algebras introduced in Section~\ref{sec:TLboundaryextns} arise as quotients of the Hecke algebras defined above. We first introduce the following notation of which we will make extensive use in the following. 
\begin{defn}
The q-number $[n]$ is defined by
\begin{equation*}
[n]=\frac{q^n -q^{-n}}{q-q^{-1}}.
\end{equation*}
\end{defn}
\noindent
We have the following proposition,
\bigskip

\begin{prop}
\label{prop:HeckeQuotients}
Let $q^{\pm 1}, q^{\pm \omega_1}, q^{\pm \omega_2}$ be such that:
\be
\label{eqn:TLparams}
\delta=[2], \quad s_1=\f{[\omega_1]}{[\omega_1+1]},
\quad s_2=\f{[\omega_2]}{[\omega_2+1]}.
\ee
Then there is a surjective algebra homomorphism, $\pi$, which is given by:
\begin{align}
\label{eqn:2BTLhomo}
\pi\left(g_i^{\pm 1}\right) &= e_i -q^{\mp 1}, \nonumber \\
\pi \left(g_0^{\pm 1}\right) &= q^{\pm \omega_1} - \left(q^{\pm(1 + \omega_1)}
- q^{\mp(1 + \omega_1)} \right) e_0, \\ 
\pi \left(g_N^{\pm 1}\right) &= q^{\pm \omega_2} - \left(q^{\pm(1 + \omega_2)}
- q^{\mp(1 + \omega_2)} \right) e_N, \nonumber 
\end{align}
\end{prop}
\begin{proof}
This follows using the 2BTL relations.
\end{proof}
\bigskip
\begin{remark}
Due to the homomorphism $\pi$, it will often be convenient for us to use the Temperley-Lieb indeterminates $\delta$, $s_1$, $s_2$ and their Hecke counterparts $q^{\pm 1}$, $q^{\pm\omega_1}$, $q^{\pm\omega_2}$ interchangeably through the identification \eqref{eqn:TLparams}.
\end{remark}

\begin{lemma}
The kernel of the homomorphism $\pi$, defined in \eqref{eqn:2BTLhomo}, is the ideal generated by the relations:
\begin{align*}
g_i g_{i+1} g_i + q^{-1} g_i g_{i+1} + q^{-1} g_{i+1} g_i + q^{-2} g_i + q^{-2} g_{i+1} + q^{-3} &=0 \qquad 1 \le i \le N-2
\end{align*}
\begin{multline}
\label{eqn:1BTLHeckeQuot}
g_1 g_0 g_1 + q^{-1} g_0 g_1 +  q^{-1} g_1 g_0 - q^{-1} \left(q^{\omega_1} + q^{-\omega_1} \right) g_1 + q^{-2} g_0 - q^{-2} \left(q^{\omega_1} + q^{-\omega_1} \right) = 0
\end{multline}
\begin{multline*}
g_{N-1} g_N g_{N-1} + q^{-1} g_N g_{N-1} +  q^{-1} g_{N-1} g_N - q^{-1} \left(q^{\omega_2} + q^{-\omega_2} \right) g_{N-1} + q^{-2} g_N \\ 
- q^{-2} \left(q^{\omega_2} + q^{-\omega_2} \right) = 0
\end{multline*}
\end{lemma}

\begin{proof}
The relations defining the 2BTL quotient of the Hecke algebra are $e_{i}e_{i\pm 1}e_i=e_i$ for $i=1,\ldots,N-1$, see Definitions~\ref{defn:TL}, \ref{defn:1BTL} and \ref{defn:2BTL}. These relations hold if and only if equations \eqref{eqn:1BTLHeckeQuot} are satisfied.
\end{proof}
\begin{remark}
\label{remark:2BTLSymmetric}
The transformation $\omega_1 \leftrightarrow -\omega_1$ is an obvious symmetry of the type B Hecke algebra given in Definition \ref{defn:HeckeB}. The Hecke quotient \eqref{eqn:1BTLHeckeQuot} is also invariant under this symmetry. Although this symmetry is broken by the definition of the boundary generator \eqref{eqn:2BTLhomo} we shall find that it 
re-emerges later -- see Remark \ref{remark:GramDetSymmetries}.
\end{remark}
The finite dimensionality of the TL and 1BTL algebras are an immediate
consequence of the finite dimensionality of the corresponding Hecke
algebras. In contrast, the 2BTL algebra is infinite dimensional. For example, for $N=2$, words of the form $\left(e_1 e_0 e_2 \right)^n$ cannot be reduced. This
will become much clearer in the diagrammatic representation which we explain in
Section~\ref{sec:DiagrammaticRepn}. 

\subsection{Integrable lattice models}

The Hecke algebras and their Temperley-Lieb quotients play an important
role in the theory of exactly solvable lattice models in statistical 
mechanics. In these so-called integrable systems the most fundamental
objects are the $R$ and $K$-matrices which arise as finite dimensional representations of
algebraic operators. The $R$ and $K$ operators defined below
satisfy the Yang-Baxter and reflection equations as a consequence of
their algebraic definition. These equations will also play a crucial role in
constructing orthogonal bases of the Temperley-Lieb algebras.
 
In addition to the identification
\[
\delta=[2], \quad s_1=\f{[\omega_1]}{[\omega_1+1]}, \quad s_2=\f{[\omega_2]}{[\omega_2+1]},
\]
in this section we will think of these indeterminates as complex numbers, and consider the 2BTL as defined over $\C$. The $R$ and $K$ operators are given in terms of the algebraic generators:
\begin{defn}
\label{defn:RK}
Let $q^u, q^\theta, q^{\bar{\theta}}\in\C$ be indeterminates. The functions $r(u)$, $k(u)$, and $\bar{k}(u)$ are defined by:
\begin{align*}
r(u) &= \frac{[u+1]}{[u]} \\
k(u) &=
-\frac{[(u-\omega_2+\theta)/2][(u-\omega_2-\theta)/2]}{[u][\omega_2+1]} \\
\bar{k}(u) &=
-\frac{[(u-\omega_1+\bar{\theta})/2][(u-\omega_1-\bar{\theta})/2]}{[u][\omega_1+1]}
\end{align*}
The R-operator is defined by
\[
R_i(u) = e_i -r(u)
\]
and the K-operators are defined by
\begin{align*}
&K_0(u)= e_0 - \bar{k}(u), &K_N(u) = e_N - k(u) \nonumber
\end{align*}
\end{defn}
\begin{prop}
\label{prop:YBE}
The $R$ and $K$-operators satisfy the Yang-Baxter \cite{B82} and left and right
reflection equations \cite{S88}:
\begin{align*}
R_i(u) R_{i+1}(u+v) R_i(v) &= R_{i+1}(v) R_i(u+v) R_{i+1}(u), \\
K_0(2v) R_1(u+v) K_0(2u) R_1(u-v)&=R_1(u-v) K_0(2u) R_1(u+v) K_0(2v), \\
K_N(2v) R_{N-1}(u+v) K_N(2u) R_{N-1}(u-v)&= R_{N-1}(u-v)K_N(2u)R_{N-1}(u+v)K_N(2v)
\end{align*}
as well as the unitarity relations:
\begin{align*}
R_i(u) R_i(-u) &= r(u)r(-u) \\
K_0(u) K_0(-u)&= \bar{k}(u)\bar{k}(-u)\\
K_N(u) K_N(-u) &= k(u)k(-u)\nonumber
\end{align*}
\end{prop}
\begin{proof}
This follows by direct application of the 2BTL relations.
\end{proof}
\begin{remark}
The quantities $K_0(u)$ and $\bar{k}(u)$ have been provided here for completeness and will play no further role in this paper. The parameter $\theta$ will
turn out to play an important role in the study of the irreducible
representations of the 2BTL, see Section~\ref{sec:irreps}.
\end{remark}

The 2BTL algebra first appeared in the study of integrable lattice
models primarily due to the existence of the following representation
\cite{deGP04,deGNPR05,N06a,N06b}:  
\begin{defn}
We define the Heisenberg spin chain to be the $2^N$ dimensional space:
\begin{align*}
\bigotimes_{i=1}^N V_i,\qquad V_i \cong \mathbb{C}^2.
\end{align*}
The term `site $i$' will be used to refer to the $i^{th}$ factor in the tensor product.
\end{defn}
\begin{prop}
\label{prop:2BTLSpinChainRepn}
The following is a family of 2BTL representations, parameterised by
$\theta$, on the Heisenberg spin chain: 
\begin{align*}
e_0 &= \f{1}{q^{\omega_1+1} - q^{-\omega_1-1}} \left\{ \sigma^+_1 - \sigma^-_1 -\half \left(q^{\omega_1}+ q^{-\omega_1} \right) \sigma^z_1 + \half \left(q^{\omega_1}-q^{-\omega_1} \right) \right\}  \\
e_i &= \sigma^+_i \sigma^-_{i+1} + \sigma^-_i \sigma^+_{i+1} + \f{q+q^{-1}}{4} \left( \sigma^z_i \sigma^z_{i+1} -1 \right) + \f{q-q^{-1}}{4} \left(\sigma^z_i - \sigma^z_{i+1} \right) \\
e_N &= \f{1}{q^{\omega_2+1} - q^{-\omega_2-1}} \left\{ - q^{\theta} \sigma^+_N + q^{-\theta} \sigma^-_N +\half \left(q^{\omega_2}+ q^{-\omega_2} \right) \sigma^z_N + \half \left(q^{\omega_2}-q^{-\omega_2} \right) \right\} 
\end{align*}
The $\sigma_i$ action is non-trivial only on site $i$ and is given by:
\begin{align*}
& \sigma^+ = \left( 
\begin{array}{cc}
0 & 1 \\
0 & 0
\end{array} \right)
&& \sigma^- = \left( 
\begin{array}{cc}
0 & 0 \\
1 & 0
\end{array} \right)
&&\sigma^z =\left( 
\begin{array}{cc}
1 & 0 \\
0 & -1
\end{array} \right)
\end{align*}
\end{prop}
\begin{proof}
This follows by direct calculation.
\end{proof}
For $i=1,2, \cdots, N-1$ the generators $e_i$ are invariant under the $U(1)$ symmetry $e_i \rightarrow U e_i U^{-1}$ with:
\[
U=\left( \begin{array}{cc}
\alpha & 0 \\
0 & \alpha^{-1}
\end{array} \right) \otimes \left( \begin{array}{cc}
\alpha & 0 \\
0 & \alpha^{-1}\end{array} \right) \otimes \cdots \otimes \left( \begin{array}{cc}
\alpha & 0 \\
0 & \alpha^{-1}\end{array} \right) 
\]
The variable $\theta$ can be thought of physically as a relative
twist of the off-diagonal terms at both ends. We shall discuss this representation further
in Section \ref{sec:SpinChainRepn}.
\section{Diagrammatic representation}
\label{sec:DiagrammaticRepn}
\subsection{The diagrammatic representation}
\begin{figure}[t]
\centerline{\includegraphics[width=0.3\textwidth]{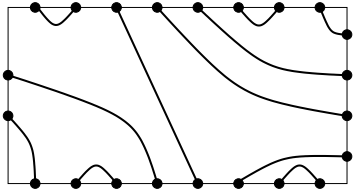}}
\caption{\label{fig:diagram1} A diagram with eight marked points on its upper and lower
  edges. The number of points on the left and right sides both have to
  be even, but can vary. The diagram above has two points on the left,
  and four on the right hand side.}
\end{figure}

\begin{figure}[t]
\begin{center}
$$
\raisebox{-60pt}{\includegraphics[height=120pt]{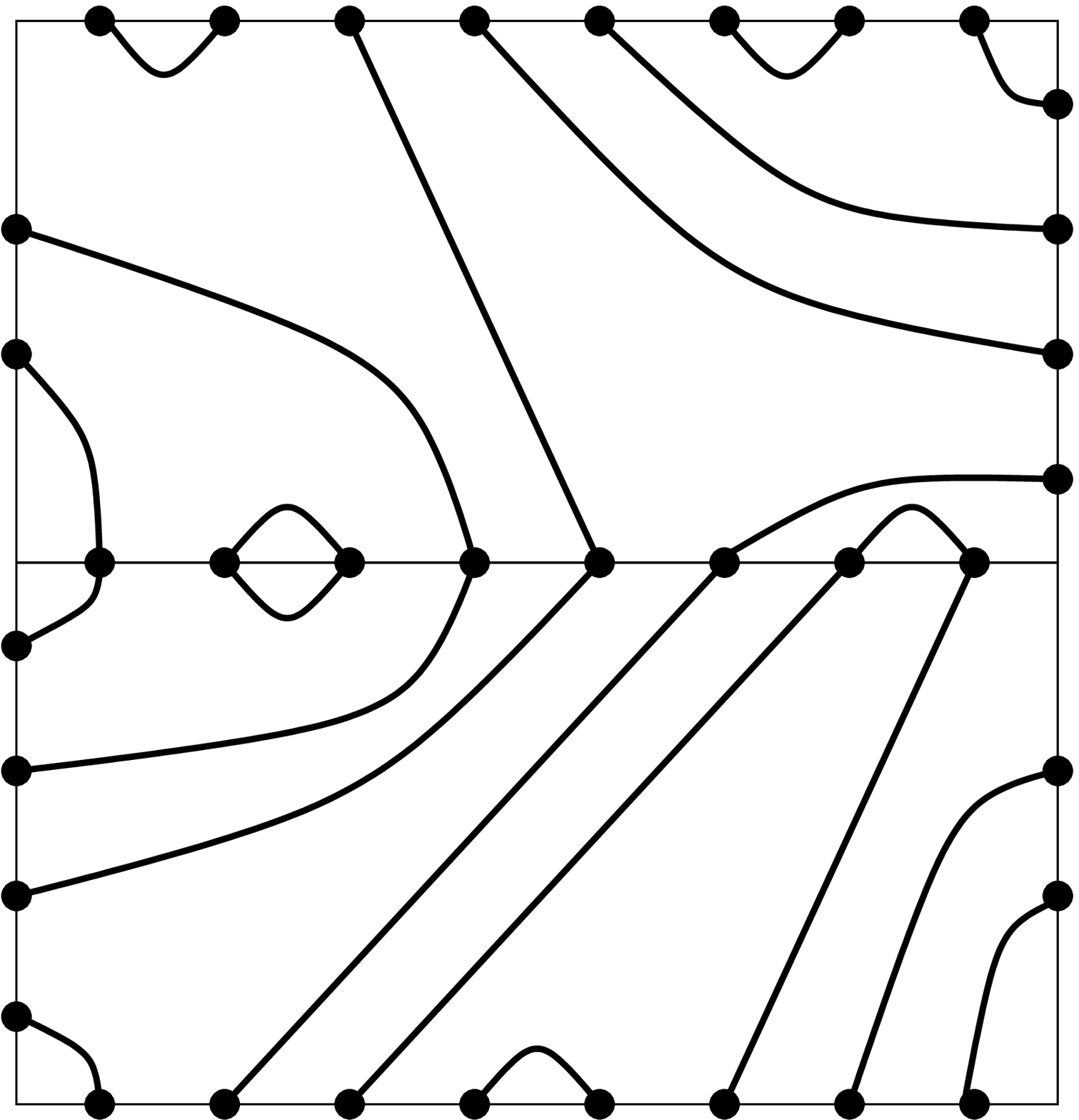}} =
\delta s_1 \ \raisebox{-30pt}{\includegraphics[height=60pt]{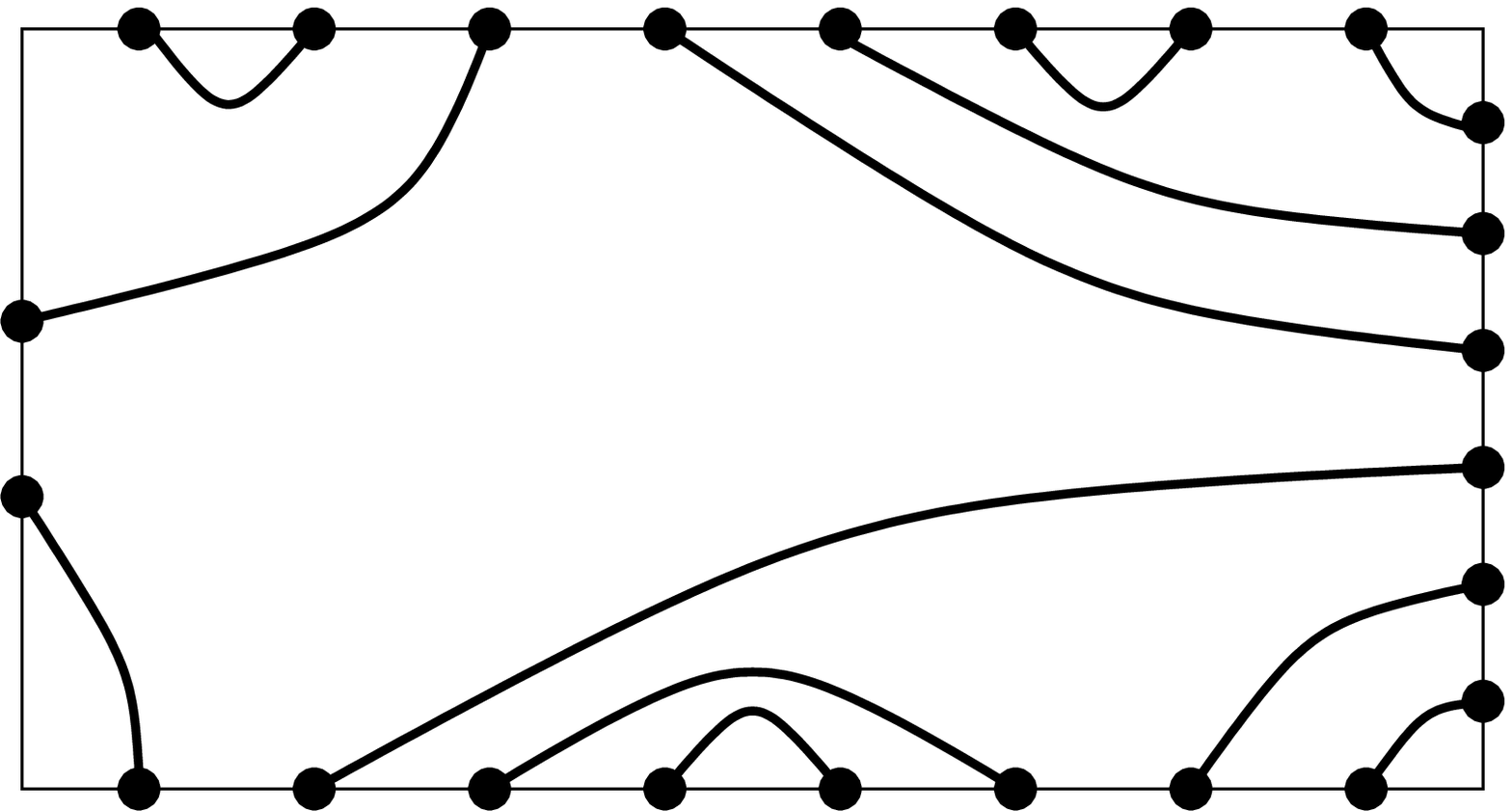}}  
$$
\end{center}
\caption{\label{fig:DiagramComposition} Composition of two diagrams consists by placing them on top
  of each other and applying the rules of Definition \ref{def:DiagramComposition}}
\end{figure}

%
%
In this section we shall give a diagrammatic representation of the
2BTL algebra. We will not attempt to include all detail in this
section, but instead refer the interested reader to \cite{MGP06}.
\begin{defn}
Take a rectangle with $N$ marked points on its upper and lower edges
and an even number of marked points on both the left and right
sides. We draw non-intersecting arcs between pairs of marked points
using each marked point once. Horizontal lines connecting the left and right side are permitted.

The set of reduced diagrams is defined
to be the subset of all such diagrams obeying the following
restrictions: 
\begin{itemize}
\item No arc has both endpoints on the left side.
\item No arc has both endpoints on the right side.
\end{itemize}
\end{defn}
An example of a reduced diagram, with $N=8$, is given in Figure \ref{fig:diagram1}. 
We consider two diagrams to be equivalent if they can be related by a smooth invertible map which preserves the edges.
\begin{defn}
\label{defn:DiagramTranspose} 
The transpose operation $^T$ operates on a reduced diagram by
reflection about the horizontal axis. 
\end{defn}
\begin{defn}
Non-reduced diagrams may contain arcs connecting the left- or right hand side to itself. For such diagrams we define an arc that connects the left (right) side to itself to be odd or even in the following way: count the number of marked points on the left (right) side below the lowest point of the arc and assign odd or even depending on its parity.
\end{defn}
As the number of marked points on each side is even we could have equivalently chosen to base the parity on the number of points above the highest point of the arc.
\begin{defn}
\label{def:DiagramComposition}
Given two reduced diagrams $A$ and $B$ the composition $A B$ is
defined by placing $A$ directly below $B$, identifying the marked
points on the common edge, and applying the following rules:
\begin{itemize}
\item Closed loops are removed with a factor $\delta$.
\item Even arcs are removed with a factor $1$.
\item Odd arcs to the left side are removed with factor $s_1$.
\item Odd arcs to the right side are removed with factor $s_2$.
\end{itemize}
\end{defn}
These rules are illustrated in Figure~\ref{fig:DiagramComposition}.
\begin{prop}
\label{diagram2alg}
The set of all reduced diagrams together with the above rules for
composition defines an algebra which is isomorphic to the 2BTL
algebra. 
\end{prop}
\begin{proof}
We identify the fundamental generators of the 2BTL algebra with the following reduced diagrams:
%
%
\begin{align*}
&& e_0 = \raisebox{-15pt}{\includegraphics[height=30pt]{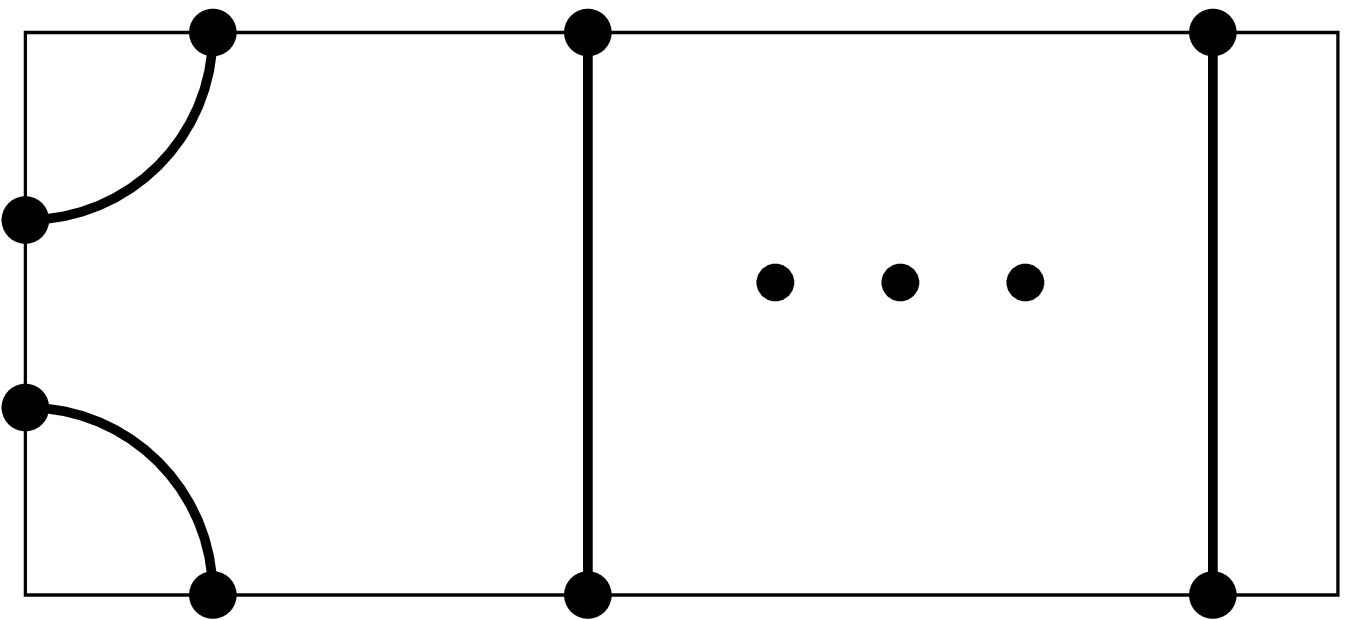}} 
&&e_i = \raisebox{-25pt}{\includegraphics[height=40pt]{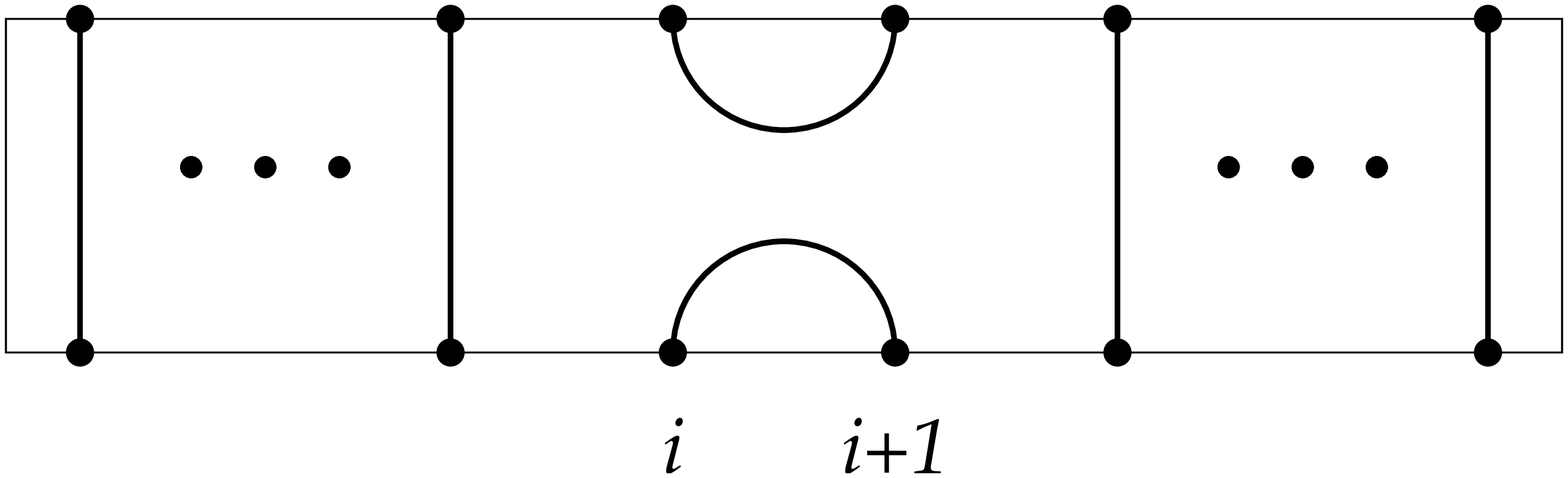}} 
&&e_N = \raisebox{-15pt}{\includegraphics[height=30pt]{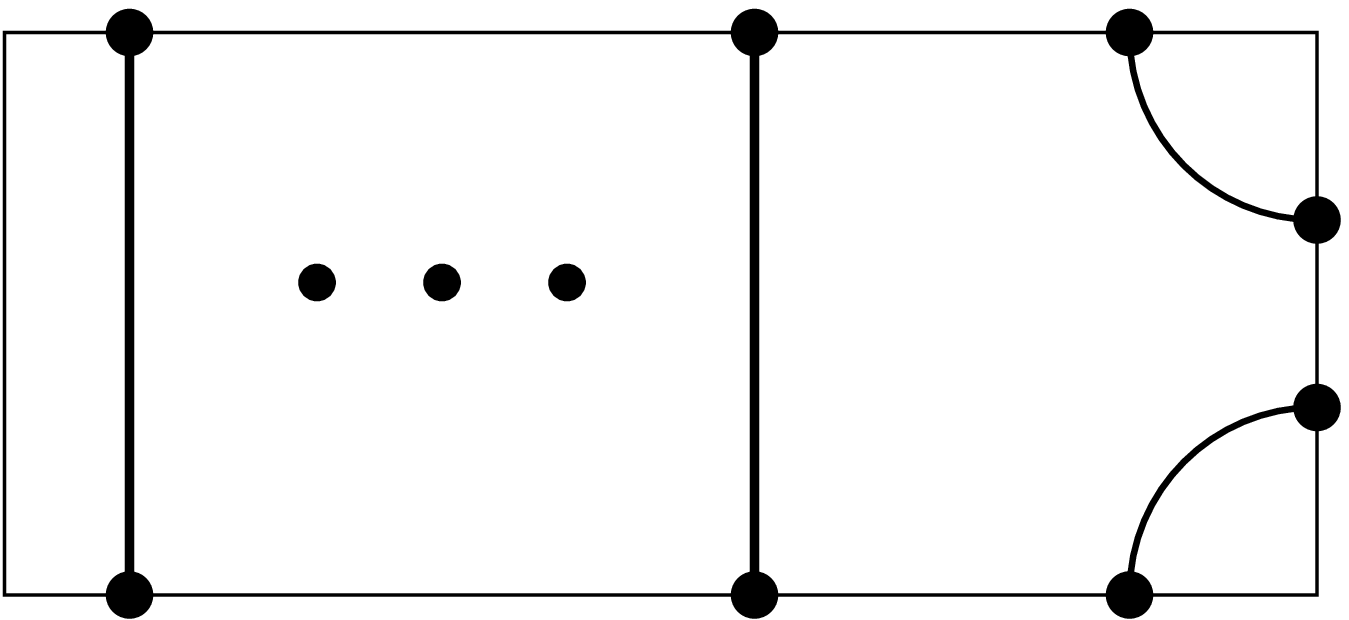}}
\end{align*}
It is easily checked that the 2BTL relations for the generators are satisfied. Furthermore, the composition rules in Definition~\ref{def:DiagramComposition} precisely correspond to word reduction in the 2BTL algebra. Hence, each reduced word in 2BTL corresponds to a reduced diagram.

Conversely, given a reduced diagram we may add closed loops that do not cross any of the other lines at the cost of scalar multiples $\delta$. In a similar way me may add even or odd arcs to either side at the cost of scalar multiples $s_1$ or $s_2$. As the number of loops ending on either side is an even number, and no loop line crosses another one, we can add closed loops and arcs in such a way that we can divide up the resulting diagram into horizontal slices where each slice is isotopic to one of the reduced diagrams corresponding to the generators $e_i$ for $i=0,\ldots,N$.
\end{proof}
\noindent
Following \cite{GL96} we might have taken the diagram algebra as an alternative definition of the 2BTL algebra. This in fact is the approach taken in \cite{MGP06}. 

\begin{figure}[h]
\begin{center}
$$
\raisebox{-60pt}{\includegraphics[height=120pt]{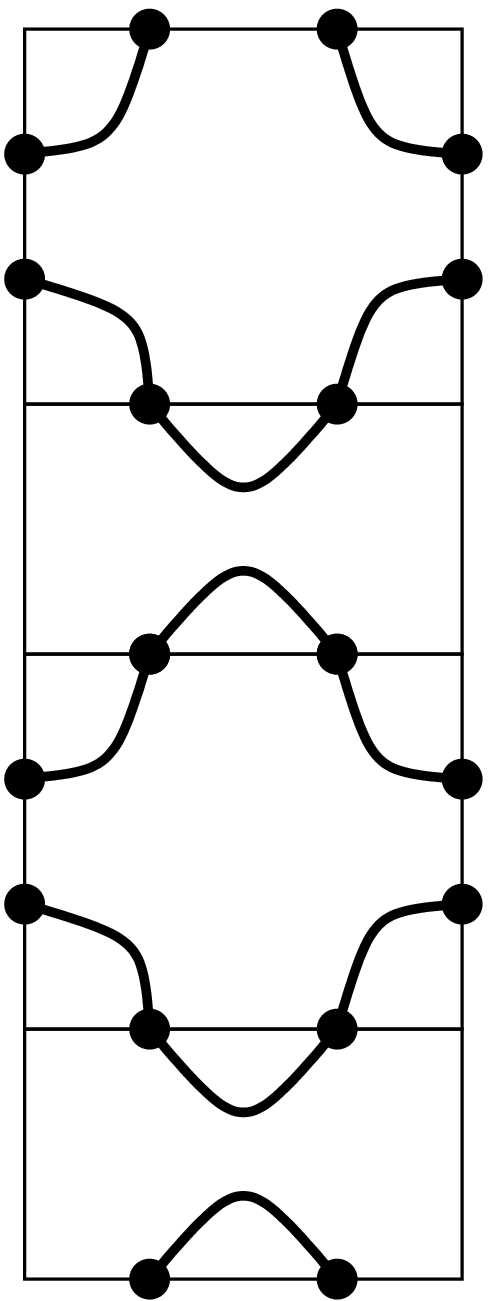}}\ = \  
\raisebox{-30pt}{\includegraphics[height=60pt]{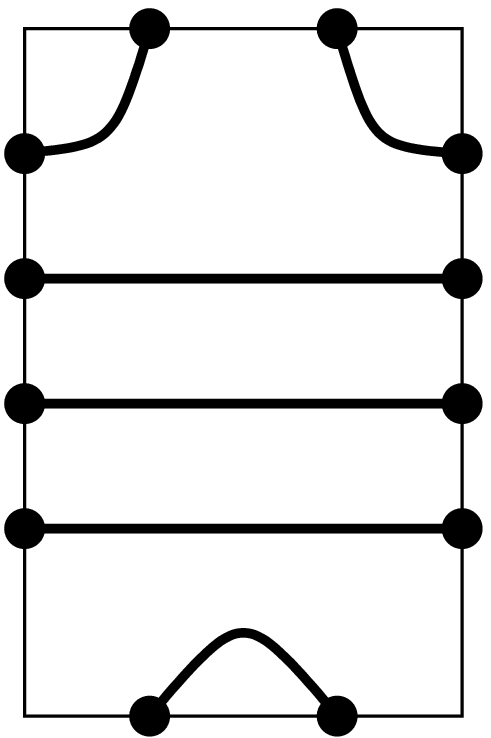}}  
$$
\end{center}
\caption{\label{fig:e1e0e2} The diagram corresponding to the word $(e_1e_0e_2)^2$ has
  three horizontal lines which cannot be removed by applying the
  algebraic rules.}
\end{figure}

%
The diagrammatic representation is infinite dimensional as there is no
restriction on the number of horizontal lines connecting the left and
right sides. For example at $N=2$ acting with $\left(e_1 e_0 e_2
\right)^n$ produces $2n-1$ horizontal lines which cannot be removed,
see Figure~\ref{fig:e1e0e2}.

\subsection{Finite dimensional quotients}
By definition reduced diagrams do not have connections from the left boundary to itself or from the right boundary to itself and therefore any horizontal lines must appear next to each other. From the reduced diagrams we form a finite dimensional subset by removing all pairs of horizontal lines. As pairs of horizontal lines are removed there always remains an even number of marked points on the left and right sides. We call this the `double quotient' of the 2BTL algebra. 

The term `double quotient' becomes apparent when this rule in expressed in terms of relations between words. We first define the idempotents $I_1$ and $I_2$ by
\begin{itemize}
\item{$N$ even}
\begin{align}
&I_1=e_1e_3\cdots e_{N-1}, && I_2=e_0 e_2 e_4\cdots e_{N-2} e_N.
\label{eqn:IdempotentsEven}
\end{align}
\item{$N$ odd}
\begin{align}
&I_1 = e_1 e_3 \cdots e_{N-2} e_N, && I_2 = e_0 e_2 \cdots e_{N-1}.
\label{eqn:IdempotentsOdd}
\end{align}
\end{itemize}
\begin{defn}
\label{defn:DoubleQuotient}
Let $b\in\C$ be an indeterminate. The double quotient of the 2BTL algebra over $\Z[\delta,s_1,s_2,b]$ is the algebra defined in Definition~\ref{defn:2BTL} subject to the additional relations:
\begin{align*}
&I_1 I_2 I_1 = b I_1, && I_2 I_1 I_2 = b I_2.
\end{align*}
\end{defn}
For example, for $N=4$ the idempotents are equal to $I_1=e_1e_3$ and
$I_2=e_0e_2e_4$, and the quotient $I_1I_2I_1=bI_1$ corresponds to 
\[
\raisebox{-45pt}{\includegraphics[height=90pt]{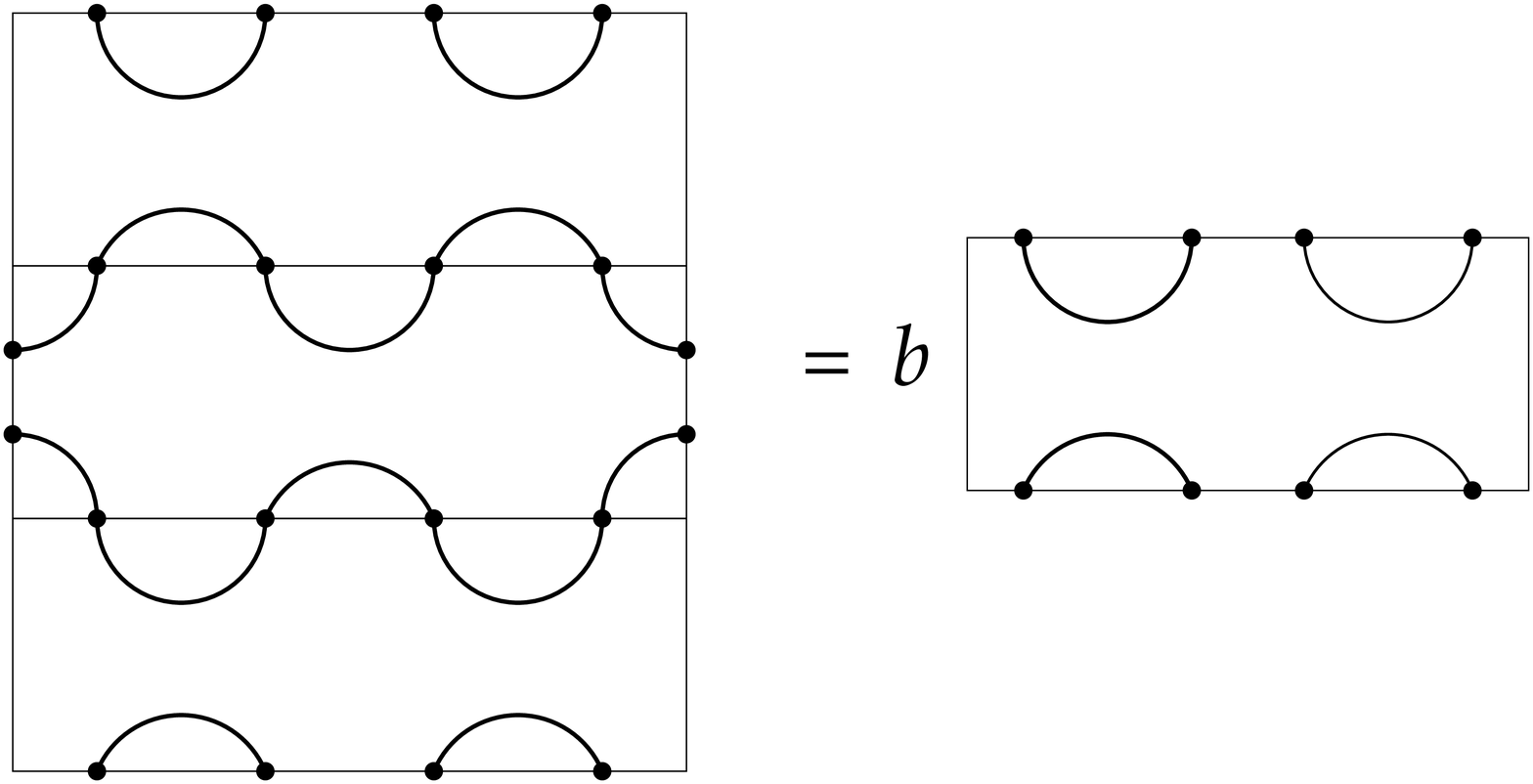}}\ .
\]
On the other hand, the quotient $I_2I_1I_2=bI_2$ is pictorially given by
\[
\raisebox{-45pt}{\includegraphics[height=90pt]{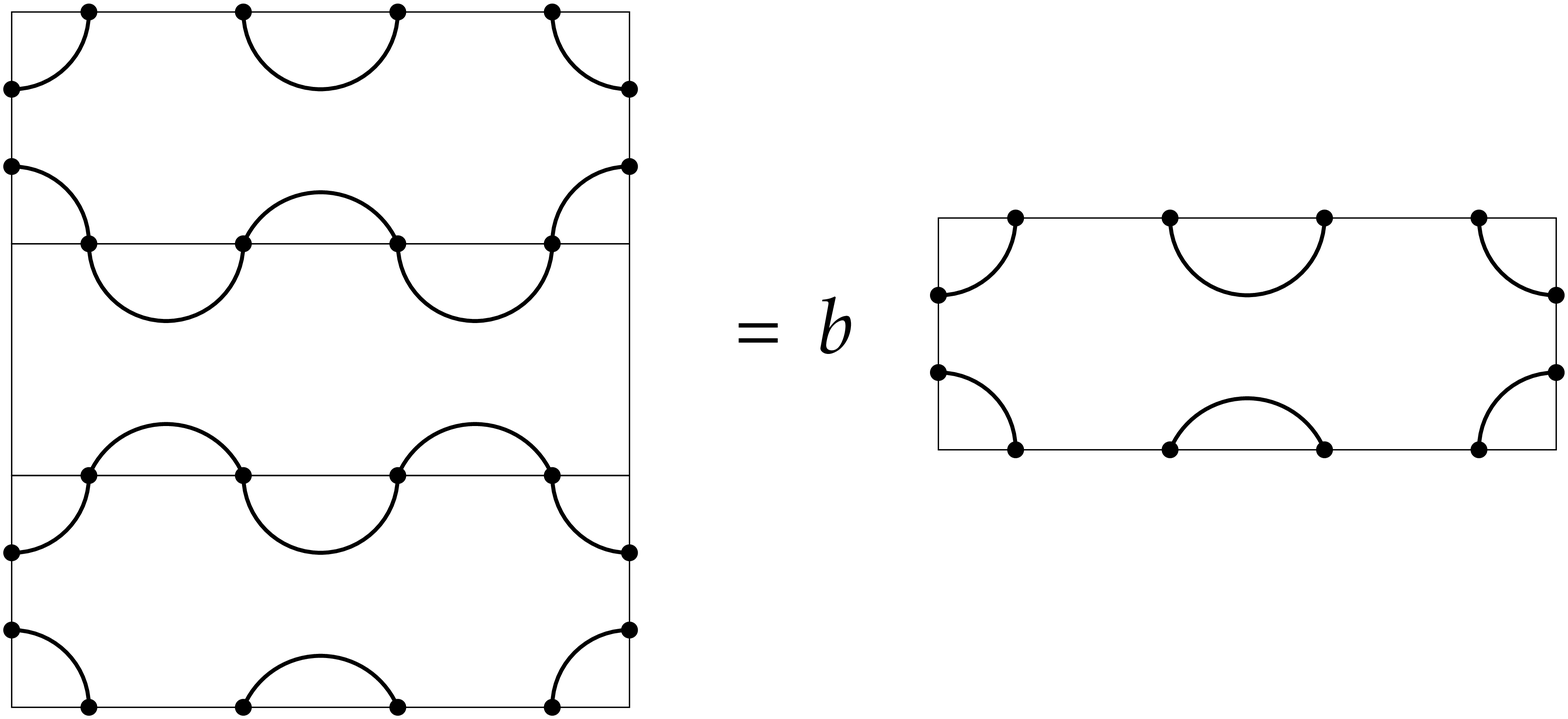}}\ .
\]
In the next section we will show that every irreducible representation
of the infinite dimensional 2BTL algebra factors through a double quotient of the form given in Definition \ref{defn:DoubleQuotient}.
\section{Irreducible representations}
\label{sec:irreps}
\subsection{Action of the centre of the 2BTL algebra}
In this paper we only study the irreducible representations of the 2BTL algebra, and we will find that all of these can be found within
finite dimensional quotients. Thinking of 2BTL as defined over $\mathbb{C}$, this is proved using the argument that, as 2BTL is countably dimensional and hence its cardinality is less than $\mathbb{C}$, Dixmier's version \cite{Dix96} of Schur's Lemma implies that all central elements are nonzero multiples of the identity in irreducible representations. In \cite{R98} this procedure was used to classify all irreducible representations of the affine Hecke algebra of type C (and all other types) for $N=2$. For our purposes, it will suffice to study the action of a particular central element of the affine Hecke algebra of type C defined by:
\be
\label{eqn:2BTLSimplestCentralElement}
Z_N = \sum_{i=0}^{N-1} \left( J^{(\rm C)}_i + \left(J^{(\rm C)}_i \right)^{-1} \right).
\ee
The main result of this section is the following:
\begin{thm}
\label{thm:IJI=bI}
Let $q^\theta\in\C$ be an indeterminate, and let $q^{\pm 1}$, $q^{\pm\omega_1}$ and $q^{\pm\omega_2}$ be such that
\[
\delta=[2], \quad s_1=\f{[\omega_1]}{[\omega_1+1]}, \quad s_2=\f{[\omega_2]}{[\omega_2+1]}.
\]
All irreducible representations of the 2BTL algebra factor through a double
quotient, given in Definition \ref{defn:DoubleQuotient}, for some specialisation of $b$. 
Moreover if the action of the central element \eqref{eqn:2BTLSimplestCentralElement} is given by
\be
Z_N =[N] \frac{[2\theta]}{[\theta]} {\bf 1},
\label{eqn:SchurZ}
\ee
then for $N$ even these irreducible representations factor through the double quotient if $\theta$ is such that
\begin{align}
b&=\frac{[(\omega_1+\omega_2+1+\theta)/2][(\omega_1+\omega_2+1-\theta)/2]}{[\omega_1+1][\omega_2+1]},
\label{eqn:beven}
\end{align}
and for $N$ odd if $\theta$ is such that 
\begin{align}
b&=-\frac{[(\omega_1-\omega_2+\theta)/2][(\omega_1-\omega_2-\theta)/2]}{[\omega_1+1][\omega_2+1]}.
\label{eqn:bodd}
\end{align}
\end{thm}
\begin{proof}
{}From the fact that $Z_N$ is central we have for any $I$ (we shall consider $I=I_1$ and $I=I_2$):
\[
I^2 Z_N = Z_N I^2 =  I Z_N I
\]
Now in the right hand side we insert the expression
\eqref{eqn:2BTLSimplestCentralElement}. In the left hand side we use \eqref{eqn:SchurZ}: 
\be
\label{eqn:CentreExpanded}
[N] \frac{[2\theta]}{[\theta]} I^2 = \sum_{i=0}^{N-1} \left( I J^{(\rm C)}_i I +  I \left(J^{(\rm C)}_i\right)^{-1} I \right) 
\ee
The main part of the proof is to simplify the expressions $I J^{(\rm C)}_i
I$. As this is not particularly illuminating we present the details in Appendix \ref{ap:IJI}.
\end{proof}
\begin{remark}
\label{remark:bNotAlwaysUnique}
The value of $b$ is uniquely determined in terms of the action of the central element
$Z$ as long as $I_1 \ne 0$ and $I_2 \ne 0$. We shall see that there exist many irreducible
representations which satisfy $I_1=0=I_2$. Such representations clearly
lie within the double quotient $I_1I_2I_1=bI_1$ and $I_2I_1I_2=bI_2$ for {\bf any}
value of $b$.
\end{remark}
\begin{remark}
Theorem \ref{thm:IJI=bI} does not imply that all irreducible representations of the 2BTL algebra can be found within the double quotient for a particular fixed value of $b$. If we were only to examine the irreducible representations for a given fixed value of $b$ we would miss a great deal of the structure of the 2BTL algebra. We shall return to this point in Section \ref{sec:OtherIrreps}.
\end{remark}
\subsection{Diagrammatic construction of irreducible representations}
\label{sec:DiagramIrreps}
\begin{defn}
A through line is an arc connecting the top and bottom edges of a diagram.
\end{defn}
The 2BTL is a cellular algebra, and to analyse its representation theory we will follow the strategy developed by Graham and Lehrer in \cite{GL96}, see also \cite{MS93,MS94}, to study its representation theory by considering the singular points of a generically nondegenerate bilinear form. In the diagrammatic representation the action of the generators either decreases the number of through lines or leaves this number unchanged. By quotienting out the former action, we can therefore construct irreducible representations on subsets of diagrams with a fixed number of through lines. Such irreducible representations are labelled by half-diagrams which we will now define.

Every full diagram, $X$, with a non-zero number of vertical through lines can
be decomposed into two half-diagrams \cite{W95,GL96}:  
\be
\label{eqn:HalfDiagramDecomp}
X = \left |   x_1 \right> \left< x_2 \right |  
\ee
where $\left |   x_1 \right>$ is the bottom part and $\left< x_2
\right |  $ is the top. We also have $X^T = \left |   x_2 \right>
\left< x_1 \right |  $. 
These half diagrams can be conveniently written using
parenthesis notation. We denote a through line by `$ | $' and
connections to the left and right by `$)$' and `$($'
respectively. For example, connected nearest neighbour marked points on the top or bottom side are denoted as `$()$'.  
\begin{example}
For $N=3$ we have:
\begin{center}
\begin{tabular}{lc}
Word & Half-diagram decomposition \\
$e_0 e_1$ &  $\Bigl|   ))  |  \Bigr > \Bigl<  ()  | \Bigr|  $ \\[10pt]
$e_0 e_3$ &  $\Bigl|   )  |  ( \Bigr > \Bigl< )  |  ( \Bigr|  $ \\[10pt]
$e_3 e_2 e_1 e_0$ &  $\Bigl  |  |(( \Bigr> \Bigr< ))| \Bigr |  $
\end{tabular}
\end{center}
\end{example}
The action of the 2BTL generators on full diagrams was given in Section
\ref{sec:DiagrammaticRepn}. For example:
\begin{center}
$e_1 \Bigl|   ))  |  \Bigr > \Bigl<  ()  | \Bigr|  = \Bigl|   ()  |  \Bigr > \Bigl<  ()  | \Bigr|$ 
\end{center}
We can now consider a vector space with basis given by the bottom parts of diagrams. For an element of this basis, say $w$, we consider a full diagram with $w$ as its bottom part. The action of the generators on $w$ is defined by considering the action of the generators on the full diagram and then extracting the bottom part using its half-diagram decomposition. For example:
\begin{center}
$e_1 \Bigl|   ))  |  \Bigr > = \Bigl|   ()  |  \Bigr > $
\end{center}
\subsubsection{Diagrams containing through lines}
If the action of a boundary generator on a diagram leaves the number
of through lines unchanged, it also preserves the parity of connections
to both the left and right hand boundary of a half-diagram. It will
therefore be convenient to define the following subsets of half-diagrams. 
\begin{defn}
\label{defn:HalfDiagramIrreps}
For a given half-diagram we define $\epsilon_1$ as the parity of connections to the left boundary and $\epsilon_2$ the corresponding quantity for the right boundary. We now define $W^{(N,n)}_{\epsilon_1,\epsilon_2}$ to be the set of all half-diagrams with parity $(\epsilon_1,\epsilon_2)$ and
$n+\frac12(\epsilon_1+\epsilon_2)$ through lines. 
\end{defn}
\begin{remark}
The action of both $I_1$ and $I_2$ on any half-diagram gives diagrams with no through lines and therefore their action vanishes on the set $W^{(N,n)}_{\epsilon_1,\epsilon_2}$ - see Remark \ref{remark:bNotAlwaysUnique}.
\end{remark}

We will see in Proposition~\ref{prop:W12irrep} that each of the sets $W^{(N,n)}_{\epsilon_1,\epsilon_2}$ when viewed as a
linear vector space forms an irreducible module of the 2BTL algebra. 
Their dimensions will be computed in what follows, and are given 
in Proposition \ref{prop:2BTLDiagramDimensions}. We shall discuss these
representations further in Section \ref{sec:OtherIrreps}.  
\begin{defn}
\label{defn:Ballot}
We define:
\[
B_{m,n} = \bin{m}{\frac{m-n}{2}} 
\]
\end{defn}
\begin{lemma}
\label{lemma:TLcounting}
The number of parenthesis sequences of length $N$ with no connections
to either boundary and $n$ through lines (these have  $N-n$ even) is
given by $B_{N,n}- B_{N,n+2}$. 
\end{lemma}
\begin{proof}
By replacing every through line `$|$' with a closing parenthesis `$)$',
we need to count all parenthesis sequences of length $N$ consisting of `$($'s and
`$)$'s with a total excess of $n$ closing parentheses, and such that at each point
the number of `$)$'s on the right is larger than the number of `$($'s. This is a
standard result which can be obtained as follows. Without the last
constraint, the total number of sequences is just $B_{N,n}$. This
overcounts the number of desired parenthesis sequences by those
sequences that at some point have one more `$($' to the right than
`$)$'s. By changing this particular `$($' into a `$)$', this set of
sequences is seen to be equinumerous to (unconstrained) parenthesis
sequences consisting of `$($'s and `$)$'s with a total excess of $n+2$
closing parentheses. Hence the Lemma follows. 
\end{proof}
\begin{lemma}
\label{lemma:1BTLcounting}
The number of parenthesis sequences of length $N$ with no connections
to the right boundary and $n$ through lines is given by
$B_{N,n}$ for $N-n$ even and $B_{N,n+1}$ for $N-n$ odd. 
\end{lemma}
\begin{proof}
In each sequence we replace every unpaired `$)$' attached to the left
boundary by a through line. This gives parenthesis sequences with no
boundary connections and extra through lines. For $N-n$ even, using
Lemma \ref{lemma:TLcounting}, we have: 
\[
\sum_{i=0}^{(N-n)/2} \left( B_{N,n+2i}- B_{N,n+2i+2} \right) = B_{N,n}
\]
The case $N-n$ odd is similar.
\end{proof}
\begin{defn}
\label{defn:2BTLIrrepSizes}
We define:
\[
M_m(n)= \sum_{i=1}^{(m+1-|n|)/2} B_{m,|n|+2i-1}
\]
\end{defn}
\begin{remark}
For $m$ odd we have $M_{m}(0)=2^{m-1}$.
\end{remark}
We are now in a position to compute the dimensions of the 2BTL modules $W^{(N,n)}_{\epsilon_1,\epsilon_2}$, whose notation now becomes apparent: 
\begin{prop}
\label{prop:2BTLDiagramDimensions}
The modules $W^{(N,n)}_{\epsilon_1,\epsilon_2}$, given in Definition 
\ref{defn:HalfDiagramIrreps}, have dimension $M_N(n)$.
%
%
\end{prop}
\begin{proof}
We consider the corresponding parenthesis sequences. We now replace
each unpaired `$($' that attaches to the right boundary with a through
line. This gives parenthesis sequences with no right boundary
connections and extra through lines. The result follows using Lemma
\ref{lemma:1BTLcounting}. 
\end{proof}
\begin{remark}
The interpretation of Lemma~\ref{lemma:1BTLcounting} is that the dimensions of the irreducible modules of the 1BTL algebra \cite{MS94,MW00} can be expressed as a sum over the dimensions of irreducible modules of the TL algebra, given in Lemma~\ref{lemma:TLcounting}. Similarly, Proposition~\ref{prop:2BTLDiagramDimensions} shows that irreducible representations of the 2BTL algebra are expressed as sums over 1BTL ones. 
\end{remark}

\subsubsection{Diagrams without through lines}
\label{sec:DiagramIrrepsNoThroughLines}
When considering the full diagrams with no through lines there is the
additional complication that some of these contain a horizontal line
(in the double quotient there cannot be more than one).

When $b\neq 0$, every full diagram can again be decomposed into two
half-diagrams, see \eqref{eqn:HalfDiagramDecomp}. These half-diagrams are obtained by
simply ignoring any horizontal lines and taking the top and bottom
parts of the full diagram as before. A horizontal line is present in
the full diagram if and only if the parities of connections to boundaries of the
upper and lower half-diagram are opposite.  

When decomposing a full diagram containing a horizontal line into half diagrams we choose to add a horizontal line to the half-diagram which contains an odd number of connections to the right
boundary (we could have equivalently chosen the left boundary). We can
then continue to use the same diagrammatic rules given in Section 
\ref{sec:DiagrammaticRepn} for the half-diagrams.

The half diagrams can once again be written using parenthesis notation
denoting connections to the left and right by `$)$' and
`$($'. For a given sequence of parentheses, the presence of a
horizontal line is completely determined. As we can have a connection to the
left or to the right at each site, it is clear that there are $2^N$
half-diagrams. We denote the corresponding space of diagrams by $W^{(N)}(b)$.
\begin{example}
For $N=3$ we have:
\begin{center}
\begin{tabular}{lc}
Word & Half-diagram decomposition \\
$e_0 e_1 e_0 e_2 e_1 e_0$ &  $\hphantom{b^{-1}}\Bigl |   ))) \Bigr> \Bigr< ))) \Bigr |  $ \\[10pt]
$e_1 e_3$ & $b^{-1} \Bigl |  \overline{()(} \Bigr> \Bigl< \overline{()(} \Bigr |  $ \\[10pt]
$e_1 e_3 e_0 e_2$ & $\hphantom{b^{-1}}\Bigl |   \overline{()(} \Bigr> \Bigl< )() \Bigr |  $ \\[10pt]
$e_0 e_2$ & $\hphantom{b^{-1}}\Bigl |   )() \Bigr> \Bigl< )() \Bigr |  $
\end{tabular}
\end{center}
\end{example}

\begin{figure}[h!]
\[
\raisebox{-90pt}{\includegraphics[height=180pt]{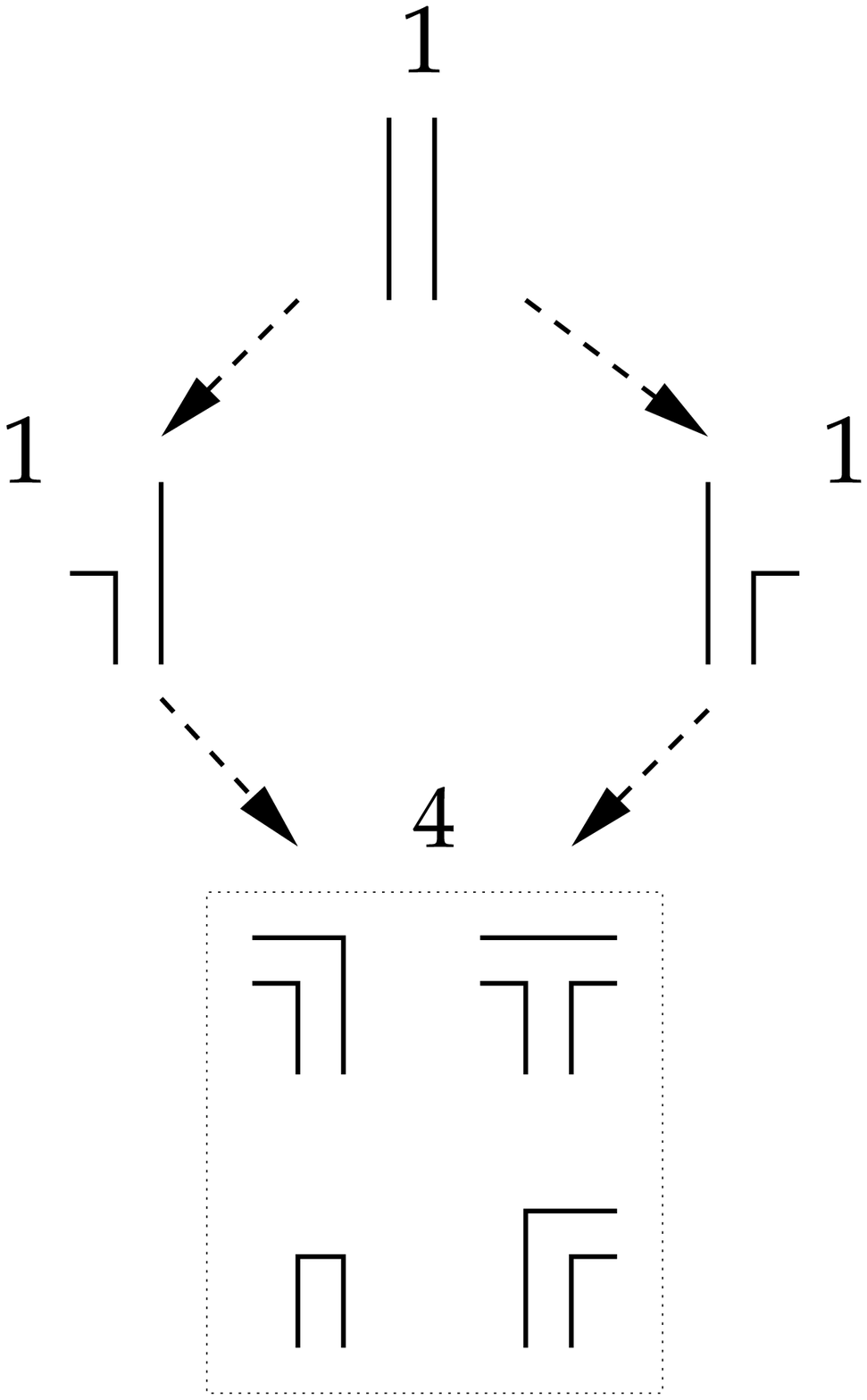}} \quad
\quad \quad \quad
\raisebox{-90pt}{\includegraphics[height=180pt]{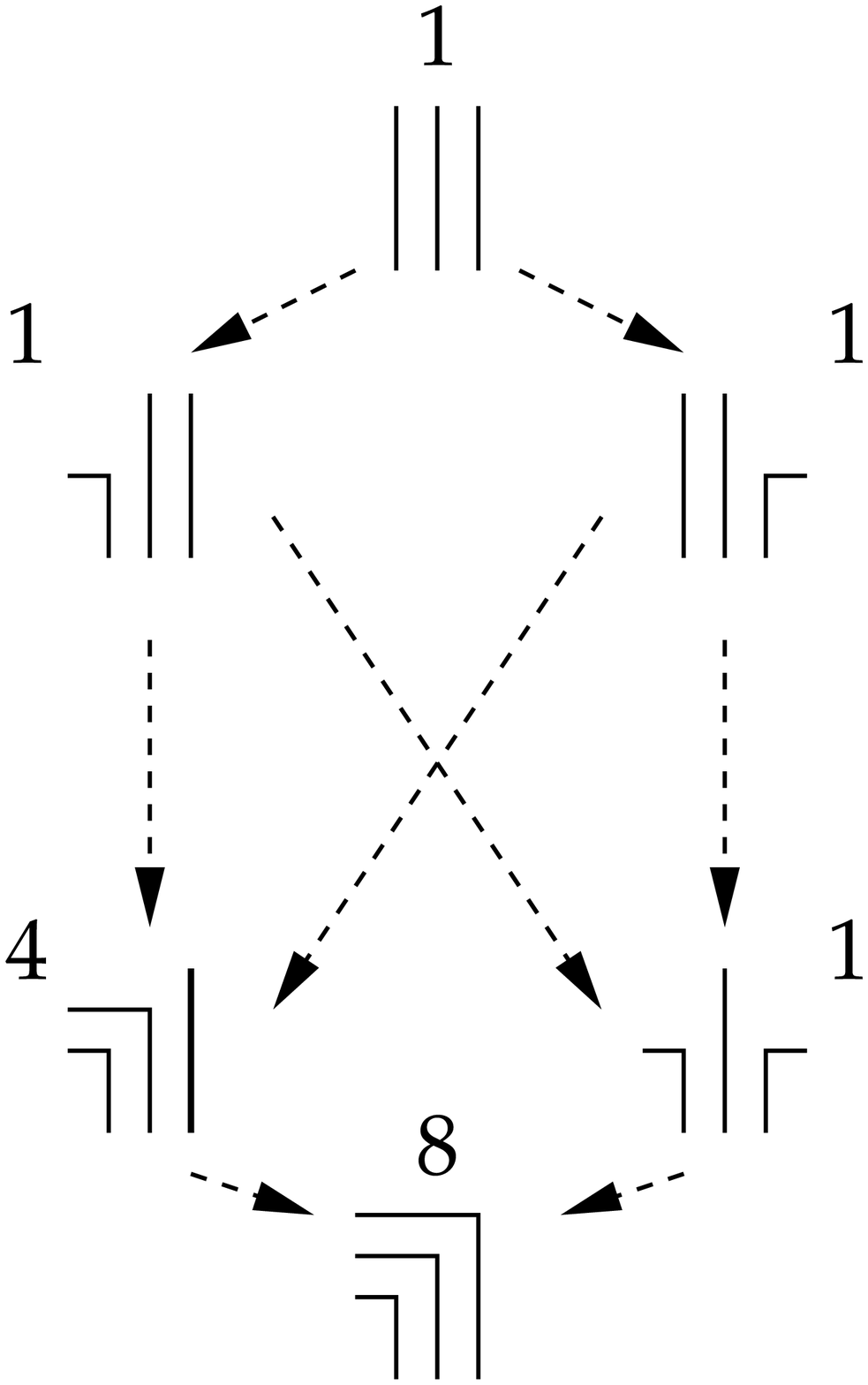}} 
\]
\caption{Organisation of the modules
  $W^{(N,n)}_{\epsilon_1,\epsilon_2}$, $W^{(N)}(b)$ and their dimensions for $N=2$ and $N=3$. For $N=2$
  all half diagrams are displayed, while for $N=3$ only representative
  diagrams for each module are depicted. The arrows indicate the
  partial ordering induced by the number of through lines.}
\label{fig:modules_N=2and3}
\end{figure}

\begin{example}
\label{example:gramN=2}
Writing shorthand $\pi$ for $| \pi \rangle$, the explicit diagrammatic representations 
for $N=2$ and $N=3$ are (see also Figure~\ref{fig:modules_N=2and3}): 

\begin{itemize}
\item For $N=2$ we have three one dimensional representations:
$W^{(2,1)}_{++}= \left\{\,|   |\,\right\}  $,
  $W^{(2,1)}_{-+}=\left\{\,) |\,\right\}  $, $W^{(2,1)}_{+-}= \left\{\,|  (\,\right\}$
and a single four dimensional representation given by:
$W^{(2)}(b)=\left \{ (),((,)),\overline{)(} \right\}$. 
\item
For $N=3$ we find four  one dimensional representations: $W^{(3,2)}_{++}= \big\{\,|||\,\big\}  $,
  $W^{(3,2)}_{-+}=\big\{\,) ||\,\big\}  $, $W^{(3,2)}_{+-}=
  \big\{\,||(\,\big\}$ and $W^{(3,2)}_{--}=
  \big\{\,)|(\,\big\}$, a four dimensional representation
  $W^{(3,0)}_{++} = \Big\{\,) )|, ()|, |(), |((\,\Big\}$ and an eight
  dimensional one $W^{(3)}(b) = \Big\{\,) )), ()), )(), )((, \overline{(((}, \overline{(()},
  \overline{()(}, \overline{))(}\,\Big\}$.
\end{itemize}
\end{example}

\begin{figure}[h!]
\[
\raisebox{-120pt}{\includegraphics[height=240pt]{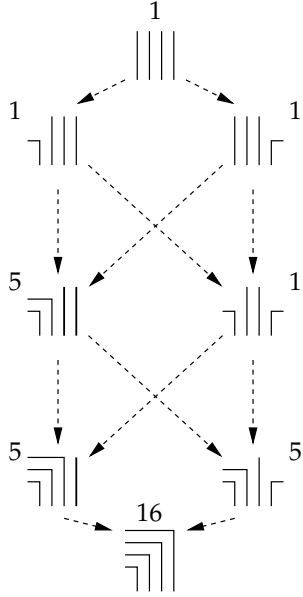}}
\]
\caption{Organisation of the modules
  $W^{(N,n)}_{\epsilon_1,\epsilon_2}$, $W^{(N)}(b)$ and their dimensions for $N=4$.}
\label{fig:module_N=4}
\end{figure}

The embedding structure of modules according to the number of
through lines can be nicely visualised. In the case of TL and 1BTL
this structure is linear, but for 2BTL we find a double linear
structure, see Figures~\ref{fig:modules_N=2and3} and
\ref{fig:module_N=4}. 

\subsection{The Gram determinant}
The representations $W_{\epsilon_1\epsilon_2}^{(N,n)}$ and $W^{(N)}(b)$ are irreducible if the set of basis vectors within each module is linearly
independent. To find when this is the case we introduce a bilinear form on the half-diagrams.
\begin{defn}
\label{defn:BilinearForm}
Similar to \cite{GL96,MS93}, we define a bilinear form $\langle\cdot |\cdot\rangle$ within the modules  $W^{(N,n)}_{\epsilon_1,\epsilon_2}$ and $W^{(N)}(b)$ in the following way:

Given two diagrams $A=|x_1\rangle\langle x_2|$ and $B=|x_3\rangle\langle x_4|$, we consider the reduced diagram $C$ obtained from their composition according to the rules in Definition \ref{def:DiagramComposition}. If $C$ is not proportional to $|x_1\rangle\langle x_4|$ we define $\langle x_2|x_3\rangle=0$. Otherwise we have:
\[
|x_1\rangle\langle x_2| \ |x_3\rangle\langle x_4| =\langle x_2|x_3\rangle\  |x_1\rangle\langle x_4|
\]
where there scalar product $\langle x_2|x_3\rangle\ $ depends only on $\langle x_2|$ and $|x_3\rangle$.
\end{defn}
\begin{example}
\label{ex:InnerProduct}
For $N=4$ and $n=2$ we have
\begin{align*}
&\Bigl\langle |()|\ \Bigr |  \ ()|| \Bigr\rangle=1,
&&\Bigl\langle ()|| \ \Bigr | \ ()|| \Bigr\rangle=\delta,
&&\Bigl\langle ||()\ \Bigr |\ ()|| \Bigr\rangle=0.
\end{align*}
\end{example}

\begin{prop}
\label{prop:W12irrep}
The modules $W_{\epsilon_1,\epsilon_2}^{N,n}$ are irreducible.
\end{prop}

\begin{proof}
We define the radical $\rad(W_{\epsilon_1,\epsilon_2}^{N,n})$ by
\[
\rad(W_{\epsilon_1,\epsilon_2}^{N,n}) = \left\{ |x\rangle\in W_{\epsilon_1,\epsilon_2}^{N,n}\;|\; \langle x | y\rangle =0\; \text{for all}\; |y\rangle\in W_{\epsilon_1,\epsilon_2}^{N,n}\right\}.
\]
Because each $|x\rangle\in W_{\epsilon_1,\epsilon_2}^{N,n}$ has the same number $n+\frac12(\epsilon_1+\epsilon_2)$ of through lines, it follows that generically $\rad(W_{\epsilon_1,\epsilon_2}^{N,n})$ is empty. As the bilinear form defined in Definition~\ref{defn:BilinearForm} is not identically zero, the Proposition follows from \cite[Prop.~3.2(ii)]{GL96}.
\end{proof}

We now define the Gram matrix in terms of this bilinear form:
\begin{defn}
\label{defn:GramMatrixDefn}
Given a basis of half-diagrams $|b_i\rangle$ for a module, the
corresponding Gram matrix $G$ is defined by $G_{ij}=\langle b_i | b_j\rangle$.
\end{defn}
It is clear from this definition that the Gram matrix is symmetric. The representations $W_{\epsilon_1\epsilon_2}^{(N,n)}$ and $W^{(N)}(b)$ will be irreducible precisely when the bilinear form is non-degenerate i.e. $\det
G\neq0$. 
\begin{example}
For the four dimensional representation $W^{(2)}(b) = \left\{ (),
)),((,\overline{)(}\right\}$ the Gram matrix $G$ reads: 
\[
G=\left(\begin{array}{cccc}
\delta & 1 & 1 & b \\
1 & s_1 & b & s_1 b \\
1 & b & s_2 & s_2 b \\
b & s_1 b & s_2 b & s_1 s_2 b
\end{array}\right)
\]
The Gram determinant is given by:
\[
\det G=b (b - s_1) (b - s_2) (b-s_1-s_2+ \delta s_1 s_2).
\]
This vanishes when $b=0$, $b=s_1$, $b=s_2$, $b=s_1+s_2-\delta s_1 s_2$. In
terms of the parameterization for $b$ given in \eqref{eqn:beven} these
points correspond to $\theta=\pm\left(1+\epsilon_1 \omega_1 + \epsilon_2 \omega_2
\right)$ where $\epsilon_1,\epsilon_2=\pm 1$. 
\end{example}

Let us now consider a particular representation with basis $\bigl\{ \rstate{b_i}\bigr \}$. Denoting both the generator and its corresponding matrix by $e_k$ we have:
\[
e_k \rstate{b_i} = \sum_{j} \left( e_k \right)_{ji} \rstate{b_j}
\]
In this paper all matrices correspond to the left action of the generators.
\begin{prop}
\label{prop:FirstGramMatrix}
On a fixed module the Gram matrix $G$ satisfies:
\begin{align*}
G e_i = e_i^T G \quad \quad i=0,1,\cdots,N 
\end{align*}
Moreover within an irreducible representation these relations are sufficient to determine the Gram matrix up to a overall scale factor.
\end{prop}
\begin{proof}
We have:
\[
\left< b_n | e_i | b_m \right> =\sum_{r} \left( e_i \right)_{rm} \left< b_n | b_r \right> = \sum_{r} \left( e_i \right)_{rm} G_{n r} = \left(G e_i\right)_{n m}.
\]
The Proposition now follows using $\left< b_n | e_i | b_m \right>=\left< b_n | e_i | b_m \right>^T$.

To prove uniqueness consider two matrices $G_1$ and $G_2$ satisfying the relations of the proposition. These must both be invertible as we are in an irreducible representation. Then we have:
\[
G_2^{-1} G_1 e_i G_1^{-1} G_2 = \left( e_i^T \right)^T = e_i
\]
and therefore $G_1^{-1} G_2$ commutes with all the $e_i$. The result now follows by Schur's Lemma as $G_1^{-1} G_2$ must be a multiple of the identity operator within any irreducible representation.
\end{proof}
In \cite{DiFGG97,DiF98} a method was developed to compute the Gram determinant for representations of the \TL algebra. The essence of this method is to give an explicit uni-triangular transformation from the word basis to a new basis in which the action of the generators is very simple and the Gram matrix is diagonal. In the next section we shall develop this approach for the 2BTL algebra. The Yang-Baxter and reflection equations together with the commutative set of Murphy elements will play a crucial role. We will obtain a basis ${\bf B}_1$ in which the action of the generators is simple, i.e. using Proposition \ref{prop:FirstGramMatrix} we will prove that the Gram matrix is diagonal in this basis.
\section{The orthogonal basis ${\bf B}_1$}
\label{sec:BasisB1}
It is our aim to compute the Gram determinant in order to understand better the irreducible representations of the 2BTL algebra. We will give conditions under which the previously described representations are irreducible and the points at which they fail to be. These points, where the Gram determinant vanishes, and the representations fail to be irreducible are called exceptional points. In this section we construct the basis which simultaneously diagonalizes all the type B Murphy elements given in
Definition \ref{defn:HeckeB}. In this basis we shall explicitly compute the action of the 2BTL
generators allowing us to show that the Gram matrix is
diagonal and hence its determinant is easily computed. 

\subsection{Construction of the basis}
We will only be concerned with the Gram matrix in the module $W^{(N)}(b)$. In order to construct the basis ${\bf B}_1$ we start with a particular idempotent $E_N$ that generates this module. The precise form of $E_N$ will be defined below. Basis elements in $W^{(N)}(b)$ are thus of the form:  
\be
\label{eqn:FundRepn}
e_{i_1} e_{i_2} \cdots e_{i_n} E_N
\ee
where $e_{i_1} e_{i_2} \cdots e_{i_n}$ is a reduced word, i.e. $n$ is the minimal number of generators $e_i$ needed to write the word \eqref{eqn:FundRepn}. We shall create a new
basis, called ${\bf B}_1$, built on $E_N$, with words of the form: 
\[
\left(e_{i_1} - \alpha_1 \right) \left(e_{i_2} - \alpha_2
\right) \cdots \left(e_{i_n} - \alpha_n \right) E_N. 
\]
A prescription for the $\alpha_i$ will be given shortly. If we order
the left ideal \eqref{eqn:FundRepn} according to the length of the
reduced words then it is clear that the change of basis to ${\bf B}_1$
is given by a lower uni-triangular matrix. They are therefore
equivalent bases and the Gram determinants computed on
either of them are equal. 
\begin{defn}
\label{defn:1BTLStartingVectors}
For $n \le N$ define $E_i$ inductively by:
\begin{align*}
E_0&= 1 \\
E_{i}&= s_1^{c_i} E_{i-1} e_{i-1} E_{i-1} \quad
\quad \quad \quad i \ge 1,
\end{align*}
where $s_1 = \f{[\omega_1]}{[\omega_1+1]}$ and $c_i=(-1)^i$.
\end{defn}
In the diagrammatic representation the $E_i$ are represented by the reduced diagrams:
%
\begin{align}
\label{eqn:StartingIdempotentDiagram}
E_i = s_1^{-\lfloor\frac{i+1}{2}\rfloor}\ \raisebox{-40pt}{\includegraphics[height=80pt]{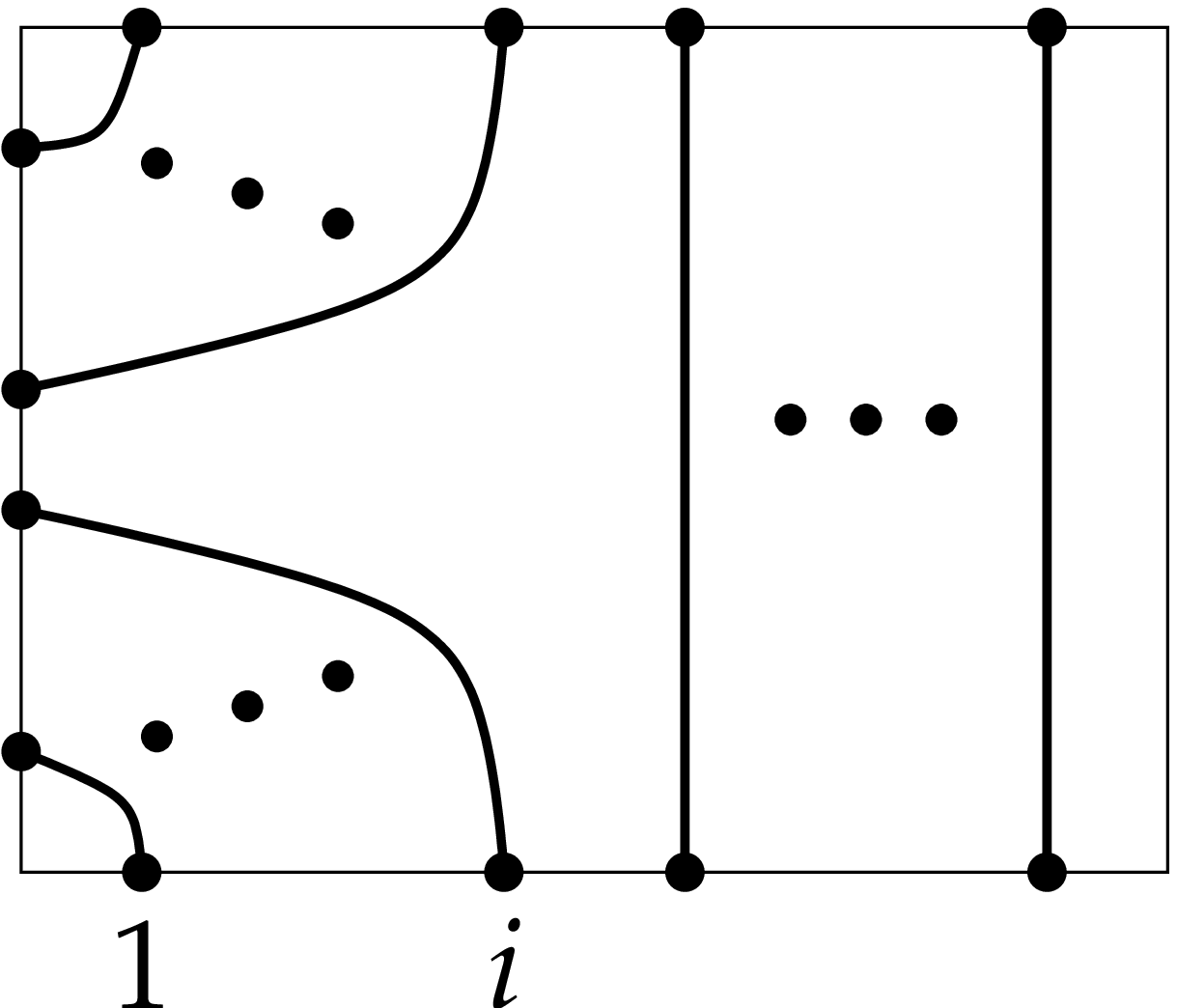}}
\end{align}
%
\noindent
It is immediate that $\left[ E_{i}, e_j
  \right]=0$ if $i<j$ and we shall use this fact repeatedly. 

\begin{figure}[t]
\centerline{\includegraphics[width=0.3\textwidth]{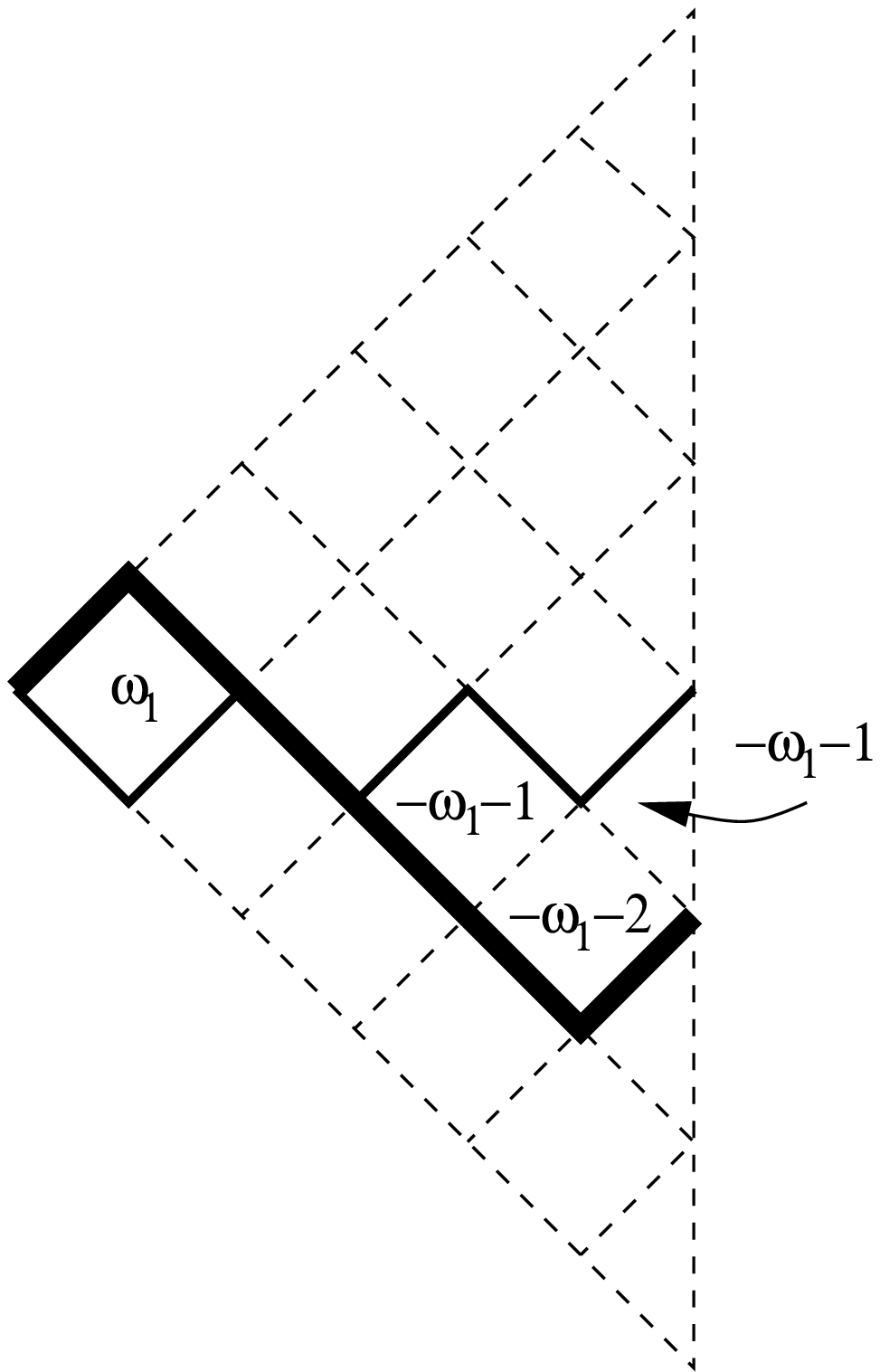}}
\caption{Paths used to label the basis elements of ${\bf B_1}$ for $N=6$. The fundamental path $p_0=(0,-1,0,-1,0,-1,0)$, which corresponds to the idempotent $E_6$, is given by the see-saw path. The thick bold path corresponds to $(0,1,0,-1,-2,-3,-2)$ and is obtained, as explained in the main text, from the fundamental path by the addition of tiles.}
\label{fig:paths}
\end{figure}

\begin{lemma} 
\label{lemma:1BTLStartingVectors}
For $j \ge i$ we have: $E_{i} E_{j}=E_{j}=E_{j} E_{i}$.
\end{lemma}
\begin{proof}
It is sufficient to prove $E_i^2=E_i$ as the remainder of the Lemma can be proved using Definition~\ref{defn:1BTLStartingVectors} and induction on $j$.

To prove $E_i^2=E_i$ we proceed by induction on $i$. The cases $i=0,1$ are easy to verify. We now assume for some fixed $i \ge 0$ that $E_i^2=E_i$ ($E_i E_{i+1}=E_{i+1}=E_{i+1}E_i$ follows from this) and $E_{i+1}^2=E_{i+1}$. Now noting that
\[ 
e_{i+1}E_{i+1}e_{i+1}= s_1^{c_{i+1}} E_i e_{i+1} E_i,
\]
we find
\[
E_{i+2}^2 = s_1^{2c_{i+2}+c_{i+1}} E_{i+1} e_{i+1} E_{i+1} = E_{i+2}
\]
where we have used Definition~\ref{defn:1BTLStartingVectors} and $c_{i+1}+c_{i+2}=0$.
\end{proof}
\begin{prop}
\label{prop:1BTLStartingVectorsMurphy}
The $E_{i}$ are eigenvectors of the type B Murphy elements $J_0^{(\rm B)}$, 
$J_1^{(\rm B)}$, $\cdots$, $J_{i-1}^{(\rm B)}$ with: 
\begin{align*}
J_{2i}^{(\rm B)} E_{2i+1} &= q^{-\omega_1 - 2i} E_{2i+1},\\
J_{2i+1}^{(\rm B)} E_{2i+2} &= q^{\omega_1 - 2i} E_{2i+2}.
\end{align*}
\end{prop}
\begin{proof}
This is a simple consequence of the following identities:
\begin{align}
\label{eqn:1BTLStartingVectorsIdentities}
J_{2i}^{(\rm B)} E_{2i} &= q^{\omega_1 - 2i} E_{2i} -
q^{-2i} \left( q^{\omega_1} - q^{-\omega_1} \right) E_{2i+1} \\ 
J_{2i+1}^{(\rm B)} E_{2i+1} &= q^{-\omega_1 - 2i-2} E_{2i+1} + q^{-2i-1}
\left( q^{\omega_1+1} - q^{-\omega_1-1} \right) E_{2i+2}, \nonumber 
\end{align}
which are proved by induction. For $i=0$ the identities in \eqref{eqn:1BTLStartingVectorsIdentities} can be explicitly checked using \eqref{eqn:2BTLhomo} and relations of the 1BTL algebra. We proceed by assuming that \eqref{eqn:1BTLStartingVectorsIdentities} holds for some fixed $i$. Multiplying the second line of \eqref{eqn:1BTLStartingVectorsIdentities} on the left and right by $g_{2i+2}$, using the recursive definition of the Murphy elements in Definition \ref{defn:HeckeB}, and finally post-multiplying by $E_{2i+2}$ we obtain: 
\begin{align}
J_{2i+2}^{(\rm B)} E_{2i+2} &= q^{-\omega_1 - 2i-2} g_{2i+2}^2 E_{2i+2} +
q^{-2i-1} \left( q^{\omega_1+1} - q^{-\omega_1-1} \right) 
g_{2i+2} E_{2i+2} g_{2i+2} E_{2i+2} \nonumber\\
&= q^{\omega_1 - 2i-2} E_{2i} - q^{-2i-2} \left( q^{\omega_1} -
q^{-\omega_1} \right) E_{2i+3}. \nonumber
\end{align}
This gives the first line of \eqref{eqn:1BTLStartingVectorsIdentities} with $i$ replaced by $i+1$.  The other case in the induction is proved in a similar fashion.
\end{proof}
We will now give a prescription of the basis ${\bf B_1}$ for the $2^N$ dimensional
irreducible representation $W^{(N)}(b)$ of the 2BTL algebra. In order to do so we will use the following convenient labelling of basis elements by paths on the tilted square lattice.
\begin{defn}
We shall write a path $p$ as a vector of local heights:
\[
p=(h_0,h_1,\cdots,h_N).
\] 
The paths are left-fixed, i.e. we have $h_0=0$, and the local heights $h_i$ are subject to the constraint $h_{i+1}-h_i=\pm 1$. We will call the path $p_0=(0,-1,0,-1,0,\cdots)$ (i.e. $h_i=0$ for
$i$ even and $h_i=-1$ for $i$ odd) the `fundamental' path.
\end{defn}
It is clear that there are precisely $2^N$ possible paths. These can all be generated recursively by adding tiles and half-tiles to the fundamental path. For example, the thick bold path in Figure~\ref{fig:paths} is obtained from the fundamental path by the addition of the three tiles and one half-tile.
%
%
%
\begin{defn}
\label{defn:BasisB1}
The first vector in the basis ${\bf B}_1=\{b_p\}$ is represented by the fundamental path $p_0$ and is given by $b_{p_0}=E_N$. The other vectors in the basis ${\bf B}_1$ are now defined recursively.

Given a path $p'$ which is obtained from another path $p$ by the addition of a tile (or half-tile) at some fixed position $i$ (i.e. $h'_i=h_i+2$), we have:
\[
b_{p'}=X_i b_{p}
\]
where the operator $X_i$ is given by the following rules:
\begin{itemize}
\item{If the tile or half-tile is added from above}
\begin{itemize}
\item Bulk: $X_i=R_i(\omega_1-h_{i-1})$
\item Right boundary: $X_i=K_N(\omega_1-h_{N-1})$
\end{itemize}
\item{If the tile or half-tile is added from below}
\begin{itemize}
\item Bulk: $X_i=R_i(-\omega_1+h_{i-1})$
\item Right boundary: $X_i=K_N(-\omega_1+h_{N-1})$
\end{itemize}
\end{itemize}
where the operators $R(u)$ and $K(u)$ were given in Definition~\ref{defn:RK}.
\end{defn}
In Figure~\ref{fig:paths} the construction of a path is illustrated in terms of tile addition. The thick bold path in this figure corresponds to the word:
\[
R_1(\omega_1) R_5(-\omega_1-2) R_4(-\omega_1-1) K_6(-\omega_1-1) E_6.
\]
In the figure we have labelled the tiles to emphasize the arguments of the operators.
%
%
%
%
\subsection{Action of the generators}
Each path $p$ corresponds to a vector $b_p \in {\bf B}_1$ and by an abuse of notation we shall often refer to the action of the generators on the paths. In the following we shall denote by $e_i \cdot p$ the action of the generator $e_i$ on the path $p$.

In order to prescribe the action of the boundary generators $e_0$ we will first need the following Lemma:
\begin{lemma}
\label{lemma:e0action}
We have:
\begin{align*}
e_0 E_N &= \f{[\omega_1]}{[\omega_1+1]} E_N, \quad \quad N> 0, \\
e_0 R_1(\omega_1) E_N  &= 0, \qquad\qquad\qquad N>1.
\end{align*}
\end{lemma}
\begin{proof}
The first of these follows trivially from the definition of $E_N$. The second is a consequence of the identity:
\[
e_0 R_1(\omega_1) e_0 e_1 = e_0  \left( e_1 - \f{[\omega_1+1]}{[\omega_1]} \right)e_0 e_1 =0
\]
\end{proof}
The behaviour of the bulk generators is very simple in the basis ${\bf B}_1$:
\begin{prop}
\label{prop:BulkActionZero}
The action of the bulk generators $e_i$ on $p\in{\bf B}_1$ vanishes when $p$ has a slope at position $i$.
\end{prop}
\begin{proof}
The proof is given in Appendix~\ref{ap:BulkActionZero}.
\end{proof}
\begin{example}
Consider the path $(0,1,0,-1,0,-1,0,\ldots)$ which has a slope at position $2$. By Definition \ref{defn:BasisB1} this corresponds to $R_1(\omega_1)E_N$. The vanishing of $e_2R_1(\omega_1)E_N$ is a consequence of 
\begin{align*}
e_2R_1(\omega_1) E_3 &= \frac{[\omega_1]^2}{[\omega_1+1]^2} e_2 \left(e_1 - \frac{[\omega_1+1]}{[\omega_1]}\right) e_0e_1e_0e_2e_1e_0=0
\end{align*}
\end{example}

\begin{thm}
\label{thm:2BTLFirstOrthogonalBasis}
The generators have the following action on paths $p=(h_0,h_1,\ldots,h_N)$ corresponding to the basis ${\bf B}_1$:
\begin{itemize}
\item Each path $p$ is an eigenvector of the left boundary generator $e_0$:
\begin{enumerate}
\item If $h_1=-1$ then $e_0 \cdot p= \f{[\omega_1]}{[\omega_1+1]} p$. 
\item If $h_1=1$ then $e_0 \cdot p= 0$.
\end{enumerate}
\item The action of the bulk generator $e_i$ is zero if we have a positive or negative slope at point $i$.

Let us take a path $p$ and another path $p'$ which is obtained from $p$ by the addition of a single tile at point $i$. The action of bulk generators $e_i$ in the basis $(p,p')$ is given by:
\begin{enumerate}
\item If $h_{i-1} \ge 0$:
\[
e_i = \left(
\begin{array}{cc}
r(\omega_1-h_{i-1}) & r(\omega_1-h_{i-1})r(-\omega_1+h_{i-1}) \\
1 & r(-\omega_1+h_{i-1})
\end{array} 
\right) 
\]
\item If $h_{i-1} < 0$:
\[
e_i = 
\left(
\begin{array}{cc}
r(-\omega_1+h_{i-1}) & r(\omega_1-h_{i-1})r(-\omega_1+h_{i-1}) \\
1 & r(\omega_1-h_{i-1})
\end{array} \right) 
\]
\end{enumerate}
\item Right boundary generator

Let us take a path $p$ and another path $p'$ which is obtained from $p$ by the addition of a right boundary half-tile. In the basis $(p,p')$ the action of $e_N$ is given by:
\begin{enumerate}
\item If $h_{N-1} \ge 0$:
\[
e_N = 
\left(
\begin{array}{cc}
k(\omega_1-h_{N-1}) & k(\omega_1-h_{N-1})k(-\omega_1+h_{N-1}) \\
1 & k(-\omega_1+h_{N-1})
\end{array} \right) 
\]
\item If $h_{N-1} < 0$:
\[
e_N = 
\left(
\begin{array}{cc}
k(-\omega_1+h_{N-1}) & k(\omega_1-h_{N-1})k(-\omega_1+h_{N-1})  \\
1 & k(\omega_1-h_{N-1})
\end{array} \right) 
\]
\end{enumerate}
\end{itemize}
\end{thm}
\begin{proof}
The action of $e_0$ follows from Lemma \ref{lemma:e0action}. By Proposition \ref{prop:BulkActionZero} the bulk generators vanish when acting on any slope. Therefore it suffices to consider their action on local maxima and local minima. The remainder of the Theorem follows directly from Definition \ref{defn:BasisB1}. For example, for $h_{i-1} \ge 0$ we have:
\begin{align*}
e_i\cdot p &= R_i(\omega_1-h_{i-1})\cdot p + r(\omega_1-h_{i-1}) p = p' + r(\omega_1-h_{i-1}) p.
\end{align*}
\end{proof}
\begin{thm}
\label{thm:TypeBMurphyEigenbasis}
The basis ${\bf B}_1$ simultaneously diagonalizes all the type B Murphy elements $J_n^{\rm(B)}$ with $n=0,1,\cdots,N-1$ and on a path $(0,h_1,h_2,\cdots,h_N)$ the eigenvalues are given by:
\[
J_{n}^{(\rm B)} = q^{\omega_1 \left( h_{n+1}-h_{n} \right)- \half \left( h_{n+1}^2-h_{n}^2 \right) +\half \left(1-2 n \right)}.
\]
\end{thm}
\begin{proof} 
The case $n=0$ follows from Theorem \ref{thm:2BTLFirstOrthogonalBasis} as the action of $J^{(\rm B)}_0=g_0$ in the basis ${\bf B}_1$ is diagonal. We now proceed by induction assuming that the Proposition is proved for all $i \le n-1$. We omit the superscript (B) for convenience.

If the path has a slope at position $n$ we have: $h_{n+1}=h_n \pm 1 = h_{n-1} \pm 2$. In this case the action of the generator $e_n$ vanishes and so $g_n=-q^{-1}$ and therefore $J_{n}=g_n J_{n-1} g_n=q^{-2} J_{n-1}$. This proves the Theorem for case $i=n$.

We now consider the cases in which $e_n$ acts non-trivially. There are two types of path to consider: $p_{\pm}=(\cdots,h_{n-1},h_{n-1}\pm 1,h_{n-1},\cdots)$ and, writing $h = h_{n-1}$, we have: 
\[
J_{n-1} \cdot p_{\pm} = q^{\pm\left(\omega_1-h\right)-n+1} p_{\pm},
\]
and need to prove that 
\[
J_{n} \cdot p_{\pm} = q^{\mp\left(\omega_1-h\right)-n+1} p_{\pm}.
\]
The ordering of these paths depends on the sign of $h$ but it is sufficient to consider $h \ge 0$. In the basis $(p_{-},p_{+})$ we have:
\[
g_n= e_n -q^{-1} =  
\left(
\begin{array}{cc}
r(\omega_1-h)-q^{-1} & r(\omega_1-h) r(-\omega_1+h) \\
1 & r(-\omega_1+h)-q^{-1}\\
\end{array} \right)
\]
and:
\[
J_{n}=g_n J_{n-1} g_n=g_n \left(
\begin{array}{cc} 
q^{-\omega_1+h+1-n} & 0 \\
0 & q^{\omega_1-h+1-n} 
\end{array} \right) g_n = \left(
\begin{array}{cc} 
q^{\omega_1-h+1-n} & 0 \\
0 & q^{-\omega_1+h+1-n}
\end{array} \right)
\]
This proves the Theorem for $i=n$.
\end{proof}
For the case of the fundamental path, corresponding to the idempotent $E_N$, Theorem \ref{thm:TypeBMurphyEigenbasis} reduces to Proposition \ref{prop:1BTLStartingVectorsMurphy}. 

{}From the action of the generators given in Theorem \ref{thm:TypeBMurphyEigenbasis} we see that the 1BTL generators act within a subspace of paths with fixed height $h_N$. We also have the following simple result:
\begin{corol}
In the basis ${\bf B}_1$ the action of the central element of the type B Hecke algebra $C_N=J^{(\rm B)}_0 J^{(\rm B)}_1 \cdots J^{(\rm B)}_{N-1}$ on a path $(0,h_1,h_2,\cdots,h_N)$ is given by:
\[
C_N = q^{h_N \omega_1 -\half h_N^2-\half N(N-2)}
\]
\end{corol}
\begin{proof}
This follows from Theorem \ref{thm:TypeBMurphyEigenbasis} using: $\sum_{i=0}^{N-1} \left(1-2i \right) = -N(N-2)$.
\end{proof}
These results suggests that the set of paths with a fixed height $h_N$ is an irreducible module for 1BTL. This is indeed the case as will be stated in Corollary~\ref{prop:1BTLcritical} of the next section.

\begin{remark}
In \cite{DJ92} the irreducible representations of the type B Hecke algebra were described using pairs of Young diagrams. The irreducible representations of the 1BTL, described here by paths of fixed height $h_N$, correspond to the restriction that both these Young diagrams have just a {\bf single} column.
\end{remark}
\subsection{Gram matrix and determinant}
\begin{prop}
In the basis ${\bf B}_1$ the Gram matrix $G$ is diagonal.
\end{prop}
\begin{proof}
It follows from Proposition \ref{prop:FirstGramMatrix} and the inductive definition of the Murphy elements given in Definition \ref{defn:HeckeB} that we have:
\[
G J_i^{(\rm B)} = \left(J_i^{(B)} \right)^T  G \quad \quad i=0,1,\cdots,N-1 
\]
In the basis ${\bf B}_1$ the Murphy elements $J_i^{(\rm B)}$ act diagonally and by Theorem \ref{thm:TypeBMurphyEigenbasis} we observe that their eigenvalues are sufficient to fully specify a path. Therefore the Gram matrix also acts diagonally on the basis ${\bf B}_1$.
\end{proof}
We now determine the eigenvalues of the Gram determinant. For convenience we define the following two functions: 
\begin{defn}
\label{defn:Functionsfg}
We define the functions $f(h)$ and $g(h)$ to be:
\begin{align*}
f(h)&=r(\omega_1-h) r(-\omega_1+h) \nonumber \\
g(h)&=k(\omega_1-h) k(-\omega_1+h) \nonumber
\end{align*}
\end{defn}
\begin{prop}
\label{prop:GramDetProcedureB1}
The eigenvalue $d_p$ of the Gram matrix for each path $p$ is given by the following recursive procedure. Let $p_0$ be the fundamental path, and let $p'$ be a path obtained from another path $p$ by the addition of a tile (or half-tile) at point $i$. The following holds:
\begin{itemize}
\item $d_{p_0}=1$.
\item If $p'$ and $p$ differ by a bulk tile we have $d_{p'} = f(h_{i-1}) d_p$.
\item If $p'$ and $p$ differ by a right boundary half-tile we have $d_{p'}=g(h_{N-1})d_p$.
\end{itemize}
\end{prop}
\begin{proof}
The fundamental path corresponds to the idempotent $E_N$ and hence it has unit norm. We now need to consider the action of the generators on the path basis. We proceed inductively by tile addition. By Theorem \ref{thm:2BTLFirstOrthogonalBasis} all the generators in the basis ${\bf B}_1$ are built of two dimensional blocks of the form:
\begin{align*}
\left(
\begin{array}{cc}
u & u v \\
1 & v \\
\end{array} \right). \nonumber
\end{align*}
It is sufficient to consider the defining relations for the Gram matrix, given in Definition \ref{defn:GramMatrixDefn}, in each of these blocks. The following identity gives recursively the entries of the Gram matrix:
\begin{align*}
\left(
\begin{array}{cc}
1 & 0 \\
0 & u v \\
\end{array} \right)
\left(
\begin{array}{cc}
u & u v \\
1  & v \\
\end{array} \right)  =  \left(
\begin{array}{cc}
u & 1  \\
u v & v \\
\end{array} \right) 
\left(
\begin{array}{cc}
1 & 0 \\
0 & u v \\
\end{array} \right)
\nonumber
\end{align*}
The specific action of the generators given in Theorem \ref{thm:2BTLFirstOrthogonalBasis} gives rise to the functions $f$ and $g$ in Definition \ref{defn:Functionsfg}.
\end{proof}
Consider again paths built by tile addition from the fundamental path $p_0$ as in Figure~\ref{fig:paths}. By Proposition \ref{prop:GramDetProcedureB1} a tile at height $h$ contributes a factor $f(h)$ and a boundary half-tile at height $h$ represents a factor $g(h)$. Let $w_p$ be the product of factors arising from all the tiles and half-tiles that must be added to the fundamental path $p_0$ to build path $p$. The determinant of the Gram matrix is given by the product of $w_p$ over all possible paths of length $N$:
\[
\det G =\prod_p w_p.
\]
\begin{lemma}
\label{lemma:PathCounting}
The number of paths of length $N$ that contain a tile at position $i$ and height $h$ is given by $2^{N-i} M_i(h)$ where $M_i(h)$ was given in Definition  \ref{defn:2BTLIrrepSizes}.
\end{lemma}
\begin{proof}
For $h \ge 0$ it is easily seen that $M_i(h)$ counts paths of length $i$ which have a  height $h+1$ or more at position $i$. The result follows as the behaviour of the path is unconstrained after position $i$. For $h < 0$ the arguments are similar but now the quantity $M_i(h)$ counts the paths of length $i$ which have height $h-1$ or less at position $i$. These are again precisely the paths that require a tile to be present at position $i$ with height $h$.
\end{proof}
\begin{thm}
\label{thm:GramDet1}
The Gram determinant is given by:
\begin{itemize}
\item for $N$ even:
\[
\det G = \alpha_N \prod_{n=0}^{(N-2)/2} \left(\prod_{\epsilon_1,\epsilon_2,\epsilon_3 =\pm 1} [(1+2n + \epsilon_1 \omega_1 + \epsilon_2 \omega_2 + \epsilon_3 \theta)/2] \right)^{M_N(2n+1)} 
\]
\item for $N$ odd:
\begin{align*}
\det G &= \alpha_N \left(\prod_{\epsilon_2,\epsilon_3=\pm 1} [(\omega_1 + \epsilon_2 \omega_2 + \epsilon_3 \theta)/2] \right)^{2^{N-1}} \\
& \quad \prod_{n=1}^{(N-1)/2} \left(\prod_{\epsilon_1,\epsilon_2,\epsilon_3=\pm 1} [(2n  + \epsilon_1 \omega_1 + \epsilon_2 \omega_2 + \epsilon_3 \theta)/2] \right)^{M_N(2n)}
\end{align*}
\end{itemize}
where $\alpha_N$ is given in both cases by:
\[
\alpha_N = \left( [\omega_1][\omega_2+1]\right)^{-2 \sum_{m=0}^{N-1} M_N(N-1-2m)}
\]
\end{thm}
\begin{proof}
The possible heights of a tile or half-tile at position $i$ is given by: $i-1-2n$ where $n=0,1, \cdots i-1$. From Lemma \ref{lemma:PathCounting} and Proposition \ref{prop:GramDetProcedureB1} we have:
\[
\det G= \left( \prod_{i=1}^{N-1} \prod_{n=0}^{i-1} f(i-1-2n)^{2^{N-i} M_i(i-1-2n)} \right) \left( \prod_{m=0}^{N-1} g(N-1-2m)^{M_N(N-1-2m)} \right)
\]
The $\theta$ independent part of this is given by:
\begin{align*}
&\prod_{i=1}^{N-1} \prod_{n=0}^{i-1} \left( \f{[\omega_1-i+2n][\omega_1-i+2n+2]}{[\omega_1-i+2n+1]^2} \right)^{2^{N-i} M_i(i-1-2n)} \\
&\quad \prod_{m=0}^{N-1} \left([\omega_1-N+2m+1] [\omega_2+1] \right)^{-2 M_N(N-2m-1)}
\end{align*}
We shall now prove that this is precisely $\alpha_N$. This is trivially true for $N=1$. We now proceed by induction on $N$ using the simple property $M_N(h)=0$ for $|h| > N$. The inductive step requires the identity:
\begin{align*}
&\prod_{n=0}^{N} \left( \f{[ \omega_1 -N +2n ]^2}{[\omega_1]^2} \right)^{M_{N+1}(N-2n)- M_N(N-1-2n) - M_N(N+1-2n)} = 1
\end{align*}
which is proved by observing that $M_{N+1}(h)=M_N(h-1)+M_N(h+1)$ for $h \ne 0$.
\end{proof}

Excluding the singular points $[\omega_i]=0$ and $[\omega_i+1]=0$, we have $\alpha_N \ne 0$ and therefore from Theorem \ref{thm:GramDet1} we deduce immediately:
\begin{corol}
\label{corol:2BTLcriticalpoints}
The $2^N$ dimensional representation $W^{N}(b)$ of the 2BTL is irreducible except at the following points:
\begin{itemize}
\item for $N$ even:
\[
\theta = \pm \left(-2n-1 + \epsilon_1 \omega_1 + \epsilon_2 \omega_2 \right) \quad n=0,1,\cdots, \f{N-2}{2} \quad \quad \epsilon_1,\epsilon_2=\pm 1 
\]
\item for $N$ odd:
\begin{align*}
\theta &= \pm \left( \omega_1 + \epsilon \omega_2 \right) \quad \quad \epsilon=\pm 1  \\
\theta &= \pm \left( -2n + \epsilon_1 \omega_1 + \epsilon_2 \omega_2 \right) \quad n=1,\cdots,\f{N-1}{2}\nonumber \quad \quad \epsilon_1,\epsilon_2=\pm 1 
\end{align*}
\end{itemize}
\end{corol}

\begin{prop}
\label{prop:1BTLcritical}
Modules defined by the set of paths of length $N$ with fixed height $h_N$ are irreducible if $\omega_1\not\in\mathbb{Z}$.
\end{prop}
\begin{proof}
We sketch the proof for $h_N\geq 0$. Reasoning along similar lines as in the proof of Theorem~\ref{thm:GramDet1}, the Gram determinant can be obtained by considering the eigenvalues associated to each tile of paths with fixed height $h_N$. Such paths lie in a tilted bounding rectangle of size $h_+\times h_-$ where $h_\pm=\frac12(N\pm h_N)$. In analogy with Lemma~\ref{lemma:PathCounting}, we now find that the number of paths of length $N$ and fixed height $h_N$, that contain a tile at position $i$ and height $h$ is given by $B_{N-i,h_N-h} M_i(h)$ where $B_{n,m}$ was given in Definition~\ref{defn:Ballot} and $M_i(h)$ in Definition~\ref{defn:2BTLIrrepSizes}.

Then, normalising such that Gram matrix eigenvalue associated to the lowest path is equal to $1$, the Gram determinant for modules defined by the set of paths with fixed height $h_N$ is given by
\[
\det G = \left( \prod_{i=1}^{N-1} \prod_{n=\max(0,i-h_+)}^{\min(i-1,h_--1)} f(h_i^*)^{B_{N-i,h_N-h_i^*} M_i(h_i^*)} \right),
\]
where we have used the abbreviation $h^*_i=i-1-2n$. The bilinear form $\langle\cdot|\cdot\rangle$ is therefore non-singular when $f(h)\neq 0$ for $h\in\mathbb{Z}$, i.e. when $\omega_1\not\in\mathbb{Z}$.
\end{proof}

\begin{remark}
\label{remark:GramDetSymmetries}
Aside from the pre-factor $\alpha_N$, the Gram determinant given in
Theorem \ref{thm:GramDet1} is invariant under $\omega_1
\leftrightarrow \omega_2$. This symmetry is manifest in the 2BTL
algebra, but we broke it in our construction of the basis ${\bf
  B}_1$. There is a further invariance under $\omega_1 \leftrightarrow - \omega_1$. As
discussed in Remark \ref{remark:2BTLSymmetric} the 1BTL quotient of
the type B Hecke algebra possesses this symmetry although it is
broken by the actual definition of the boundary generator $e_0$. 

The additional symmetry $\omega_2 \leftrightarrow \theta$ is unexpected,
and the origin of this extra symmetry will be discussed elsewhere. 
\end{remark}
\subsection{The spin chain representation}
\label{sec:SpinChainRepn}
In this subsection we shall discuss further the spin chain, or tensor
product, representation of the 2BTL algebra on the space $\bigotimes_{i=1}^N V_i$ with
$V_i\approx\mathbb{C}^2_i$, which was given in Proposition
\ref{prop:2BTLSpinChainRepn}. We prove in Theorem
\ref{thm:SpinChainIdentified} that this representation is 
equivalent to $W_N(b)$.

Using the notation
\[ \uparrow = \left(\begin{array}{@{}c@{}} 1\\0\end{array}\right),\qquad
\downarrow = \left(\begin{array}{@{}c@{}} 0\\1\end{array}\right),
\]
we make the following definition.
\begin{defn}
We define $\bar{E}_N$ inductively through:
\begin{align*}
\bar{E}_1 &= \left( q^{-\omega_1} \uparrow + \downarrow \right)  \\
\bar{E}_{2i} &= \bar{E}_{2i-1} \otimes \left( q^{\omega_1+1} \uparrow + \downarrow \right) \quad i \ge 1 \\
\bar{E}_{2i+1} &= \bar{E}_{2i} \otimes \left( q^{-\omega_1} \uparrow + \downarrow \right) \qquad i \ge 1
\end{align*}
\end{defn}
\begin{prop}
\label{prop:EEbarIdentity}
In the spin chain representation given in Proposition \ref{prop:2BTLSpinChainRepn}, and $E_i$ as in Definition \ref{defn:1BTLStartingVectors}, we have:
\[
E_i \bar{E}_N = \bar{E}_N
\]
\end{prop}
\begin{proof}
The case $i=0$ is trivial and $i=1$ is simple to prove using the
explicit action of  $e_0$ given in Proposition
\ref{prop:2BTLSpinChainRepn}. We now proceed by induction using the
definition of $E_{i}$ and the action of $e_i$, given in Proposition
\ref{prop:2BTLSpinChainRepn}, on sites $i$ and $i+1$: 
\begin{align*}
&e_i \uparrow \uparrow =0,
&&e_i \uparrow \downarrow = -q^{-1} \uparrow \downarrow + \downarrow \uparrow, 
&&e_i \downarrow \uparrow =  \uparrow \downarrow - q \downarrow \uparrow, 
&&e_i \downarrow \downarrow =0
\end{align*}
\end{proof}
The use of the variable $\theta$ in Proposition \ref{prop:2BTLSpinChainRepn} is justified by the following Lemma:
\begin{lemma}
\label{lemma:SpinSimpleBulkActionZero2}
In the spin chain representation Lemma \ref{lemma:SimpleBulkActionZero2} holds with  $E_N$ replaced by $\bar{E}_N$.
\end{lemma}
\begin{proof}
As discussed in the proof of Lemma \ref{lemma:SimpleBulkActionZero2} it is sufficient to verify:
\begin{itemize}
\item For $N$ even: $e_{N-1} K_N(-\omega_1-1) \bar{E}_N=0$
\item For $N$ odd: $e_{N-1} K_N(\omega_1) \bar{E}_N =0$.
\end{itemize}
The action of the generators $e_{N-1}$ and $e_N$, given in Proposition \ref{prop:2BTLSpinChainRepn}, is non-trivial only on the final two sites of the spin chain. The Lemma follows by direct calculation.
\end{proof}
\begin{thm}
\label{thm:SpinChainIdentified}
The spin chain representation of 2BTL, given in Proposition
\ref{prop:2BTLSpinChainRepn}, is equivalent to the $2^N$ dimensional representation $W_N(b)$.
\end{thm}
\begin{proof}
The results given in Proposition \ref{prop:1BTLStartingVectorsMurphy}, Lemma \ref{lemma:SimpleBulkActionZero1}, and Lemma \ref{lemma:e0action} were proved using only the 2BTL algebra. Therefore by acting on $\bar{E}_N$, and using Proposition \ref{prop:2BTLSpinChainRepn}, we find corresponding results are true when $E_N$ is replaced by $\bar{E}_N$. Lemma \ref{lemma:SpinSimpleBulkActionZero2} gives the result corresponding to Lemma \ref{lemma:SimpleBulkActionZero2} with $E_N$ replaced by ${\bar E}_N$.

We have therefore established that all of the fundamental results required to prove Theorem \ref{thm:2BTLFirstOrthogonalBasis} and Theorem \ref{thm:TypeBMurphyEigenbasis} hold equally well with $E_N$ replaced by $\bar{E}_N$. Therefore we can construct a new basis ${\bf \bar{B}}_1$, in a similar way to basis ${\bf B}_1$, with $\bar{E}_N$ rather than $E_N$. The action of all generators in bases ${\bf B}_1$ and ${\bf \bar{B}}_1$ is identical and therefore the representations are equivalent.
\end{proof}
\begin{remark}
The quantity $\bar{E}_N$ is part of the so-called ${\bf Q}$ basis found in a study of the 1BTL algebra \cite{NRdeG05}.
\end{remark}
\section{Other irreducible representations.}
\label{sec:OtherIrreps}
In our construction of basis ${\bf B}_1$ in Section \ref{sec:BasisB1} we only considered the $2^N$ dimensional representation of the 2BTL algebra. This representation is parameterized by a single additional number $\theta$ and is generically irreducible - see Corollary \ref{corol:2BTLcriticalpoints}.

We now consider the case where $\theta$ is not generic but takes one of the exceptional values given in Corollary \ref{corol:2BTLcriticalpoints}. We shall show that at these points there is an invariant subspace and one may obtain smaller irreducible representations of the 2BTL algebra. This can also be understood directly in the spin chain representation of Section \ref{sec:SpinChainRepn} \cite{N06b}. 
\begin{remark}
We shall assume in the following that the exceptional points given in Corollary \ref{corol:2BTLcriticalpoints} are distinct. 
\end{remark}
The discussion becomes more involved in the cases in which two or more of these exceptional points coincide. 
We shall use the following notation:
\begin{defn}
Consider a space $V$ and a linear transformation $T: V \rightarrow V$. If $T$ has a non-trivial invariant subspace $Y$ then we shall denote the appearance of this by the notation $V \longrightarrow Y$.
\end{defn}
The following Proposition gives a more explicit understanding of the points where the Gram determinant vanishes:
\begin{prop}
\label{prop:IndecStruct}
Let $P$ denote the set of paths in basis ${\bf B}_1$. Then the action of the 2BTL generators is given by:
\begin{itemize}
\item For $\theta=\pm \left(-m+\omega_1 \pm \omega_2 \right)$ and $m \ge 0$:
 $P \longrightarrow P(h_{N} \ge m+1)$.
\item For $\theta=\pm \left(-m-\omega_1 \pm \omega_2 \right)$ and $m > 0$:
 $P \longrightarrow P(h_{N} \le -m-1)$.
\end{itemize}
\end{prop}
\begin{proof}
For $\theta=\pm \left(-m+\omega_1 + \omega_2 \right)$ and $m \ge 0$ we have $k(-\omega_1+m)=0$. Now consider the two sets of paths $P_{\pm} = (0;h_1;\cdots;m ;m \pm 1)$ in basis ${\bf B}_1$. From Theorem \ref{thm:2BTLFirstOrthogonalBasis} in the basis $\left(P_+,P_-\right)$ we have:
\[
e_N = \left(
\begin{array}{cc}
\f{[\omega_2]}{[\omega_2+1]} & 0 \\
1 & 0
\end{array}
\right)
\]
and so the paths in set $P_+$ are invariant under $e_N$. By the action of the 2BTL generators on the set $P_+$ we generate all possible paths with $h_N \ge m+1$.

The other cases are treated in a similar manner.
\end{proof}
At every point where we have a reducible but indecomposable representation $V \longrightarrow Y$ we can obtain two different irreducible representations by taking $Y$ and $V/Y$.
\begin{thm}
\label{thm:2BTLIrreps}
We have the following non-equivalent irreducible representations of the 2BTL algebra:
\begin{itemize}
\item For $N$ even we have representations $V^{(N,2i+1)}_{\epsilon_1,\epsilon_2}$ and $\tilde{V}^{(N,2i+1)}_{\epsilon_1,\epsilon_2}$ where $\epsilon_1,\epsilon_2=\pm 1$ and $i=0,1,\cdots, \f{N-2}{2}$.

\item For $N$ odd we have representations $V^{(N,2i)}_{\epsilon_1,\epsilon_2}$ and $\tilde{V}^{(N,2i)}_{\epsilon_1,\epsilon_2}$ where $\epsilon_1,\epsilon_2=\pm 1$ and $i=1,\cdots, \f{N-1}{2}$ and also $V^{(N)}_{\epsilon}$ and $\tilde{V}^{(N)}_{\epsilon}$ where $\epsilon=\pm 1$.
\end{itemize}
Their dimensions are given by:
\begin{align*}
\dim V^{(N,n)}_{\epsilon_1,\epsilon_2} &= M_N(n) \\
\dim \tilde{V}^{(N,n)}_{\epsilon_1,\epsilon_2} &= 2^N- M_N(n) \\
\dim V^{(N)}_{\epsilon} &= 2^{N-1} = \dim \tilde{V}^{(N)}_{\epsilon}
\end{align*}
where $M_N(n)$ was given in Definition \ref{defn:2BTLIrrepSizes}.

The action of the centre $Z_N$, defined in \eqref{eqn:2BTLSimplestCentralElement}, is given by:
\begin{align*}
&Z_N V^{(N,n)}_{\epsilon_1,\epsilon_2}= \left[N \right] \f{[2 \left( -n +
      \epsilon_1\omega_1+\epsilon_2 \omega_2\right) ]}{[-n+\epsilon_1
      \omega_1+\epsilon_2 \omega_2]}
  V^{(N,n)}_{\epsilon_1,\epsilon_2}, &&Z_N V^{(N)}_{\epsilon} = \left[N \right] \f{[2
      \left( \omega_1+\epsilon \omega_2\right) ]}{[\omega_1+\epsilon
      \omega_2]}V^{(N)}_{\epsilon} 
\end{align*}
The action of the centre on $\tilde{V}^{(N,n)}_{\epsilon_1,\epsilon_2}$ is the same as on $V^{(N,n)}_{\epsilon_1,\epsilon_2}$ and on $\tilde{V}^{(N)}_{\epsilon}$ is the same as on $V^{(N)}_{\epsilon}$.
\end{thm}
\begin{proof}
These representations are all obtained from the $2^N$ dimensional one at the exceptional points given by Corollary \ref{corol:2BTLcriticalpoints}. Their dimensions follow immediately from Proposition \ref{prop:IndecStruct}.
\end{proof}
The fact that the dimensions of $V^{(N,n)}_{\epsilon_1,\epsilon_2}$ coincide with the dimensions of the irreducible sets of half-diagrams $W^{(N,n)}_{\epsilon_1,\epsilon_2}$ found in Proposition \ref{prop:2BTLDiagramDimensions} motivates the conjecture that all the irreducible representations found in the diagrammatic approach of Section \ref{sec:DiagramIrreps} can be described in this way:
\begin{conj}
\label{conj:DiagramIrrepsIdentified}
Let $\theta=-n+\epsilon_1\omega_1+\epsilon_2\omega_2$. Then for $n>0$ the irreducible representation $W^{(N,n)}_{\epsilon_1,\epsilon_2}$ of 2BTL is equivalent to $V^{(N,n)}_{\epsilon_1,\epsilon_2}$ and $W^{(N,0)}_{+,+}$ is equivalent to $V^{(N)}_{+}$.
\end{conj}
\begin{example}
Let $N=2$, and consider the case $n=1$ and $\epsilon_1=-\epsilon_2=+$, which corresponds to $b=s_1$. In this case the representation $W^{(2)}(b=s_1)$ has a one-dimensional invariant subspace $V_{+-}^{(2,1)}= \{\ )(-s_1((\ \}$. On this module it is easily checked that we have the one-dimensional representation given by
\[
e_0\mapsto 0, \qquad e_1\mapsto 0, \qquad e_2\mapsto s_2.
\]
It is quickly checked that this representation is exactly the same as $W_{+-}^{(2,1)} = \{\ |(\ \}$.
\end{example}
\appendix
\section{Proof of Theorem~\ref{thm:IJI=bI}}
\label{ap:IJI}
In this appendix we complete the proof of Theorem~\ref{thm:IJI=bI} by
simplifying further the expressions $I J^{(\rm C)}_i I$ in \eqref{eqn:CentreExpanded}.

We shall consider the cases of $N$ even and $N$ odd separately. In
each case the proof consists of two parts. First we prove recursion
relations and second we deal explicitly with the remaining terms. 
\subsection{Preliminary Lemmas}
\label{sec:IJILemmas}
The following identities follow from
the definitions of Section \ref{sec:TLboundaryextns}:
\begin{align}
\label{eqn:eg_relations}
e_i g_{i \pm 1} g_i &= -q^{-1} e_i e_{i \pm 1}, \qquad 
g_i g_{i \pm 1} e_i = -q^{-1} e_{i \pm 1} e_i\nonumber\\
e_N g_{N-1}^{-1} g_N g_{N-1} e_N &= e_N \left( q^{-\omega_2} +
\left( q^{\omega_2-1}-q^{-\omega_2+1} \right)e_{N-1}\right)e_N,\\
e_0g_1g_0g_1e_0 &= q^{-1} e_0 \left( q^{-1-\omega_1} + \left( q^{\omega_1} -
q^{-\omega_1} \right) e_1 \right)e_0. \nonumber
\end{align}
Similar identities hold for the inverse Hecke generators,
and differ only in interchanging $q \leftrightarrow q^{-1}$. These
identities together with the commutation relations of the Murphy
elements $J^{(\rm C)}_i$ and the generators $g_j$, will be used
repeatedly in this section without comment.  
\begin{lemma}
\label{lemma:NEvenRecursion}
For $N$ even we have:
\begin{align}
&I_1 J^{(\rm C)}_{2i+1} I_1 = q^{2} I_1 J^{(\rm C)}_{2i} I_1,
&&I_1 J^{(\rm C)}_{2i+2} I_1 = q^{-2} I_1 J^{(\rm C)}_{2i} I_1.
\label{eqn:recurJa}
\end{align}
The first equality holds for $0 \le i \le (N-2)/2$ and the second
for $0 \le i \le (N-4)/2$. We also have:
\begin{align}
&I_2 J^{(\rm C)}_{2i+2} I_2 = q^2 I_2 J^{(\rm C)}_{2i+1} I_2,
&&I_2 J^{(\rm C)}_{2i+3} I_2 = q^{-2} I_2 J^{(\rm C)}_{2i+1} I_2. 
\label{eqn:recurJb}
\end{align}
Here, the first equality holds for $0 \le i \le (N-4)/2$ and the
second for $0 \le i \le (N-6)/2$.
\end{lemma}
\begin{proof}
We use the recursive definition of $J_i$:
\[
I_1 J_{2i+1} I_1 = I_1 g_{2i+1} J_{2i} g_{2i+1} I_1 = q^{2} I_1 J_{2i} I_1,
\]
proving the first relation of \eqref{eqn:recurJa}. The second relation of \eqref{eqn:recurJa} follows from:
\[
I_1 J_{2i+2} I_1 = I_1 g_{2i+2} g_{2i+1} J_{2i} g_{2i+1} g_{2i+2} I_1 
= q^{-2} I_1 e_{2i+2} J_{2i} e_{2i+2} I_1 = q^{-2} I_1 J_{2i}I_1.
\]
The relations in \eqref{eqn:recurJb} are proved in a similar manner. 
\end{proof}
\begin{lemma}
\label{lemma:IJ0I}
For $N$ even the following identities hold: 
\begin{align*} 
I_1 J^{(\rm C)}_0 I_1 &= q^{-2} [2]^{(N-2)/2} 
\left( \frac{[2 \left(\omega_1+\omega_2+1 \right)]}{[\omega_1+\omega_2+1]} I_1
- (q-q^{-1})^2 [\omega_1+1][\omega_2+1] I_1 I_2 I_1 \right)\\  
I_2 J^{(\rm C)}_0 I_2 &= q^{-\omega_1} [2]^{(N-2)/2}\left(
q^{-\omega_2} \f{[\omega_1][\omega_2]}{[\omega_1+1][\omega_2+1]} I_2 + (q-q^{-1})[\omega_2-1] I_2 I_1 I_2 \right)\\
I_2 J^{(\rm C)}_1 I_2 &= q^{-3}  [2]^{(N-4)/2}\left( \frac{[\omega_1][\omega_2][2\left(\omega_1+\omega_2\right)]}{[\omega_1+1][\omega_2+1][\omega_1+\omega_2]} I_2 - (q - q^{-1})^2 [\omega_1] [\omega_2-1] I_2 I_1 I_2 \right) \\ 
I_2 J^{(\rm C)}_{N-1} I_2 &= q^{-\omega_2-N+1} [2]^{(N-2)/2} \left( q^{-1-\omega_1} \f{[\omega_1][\omega_2]}{[\omega_1+1][\omega_2+1]} I_2 + (q-q^{-1})[\omega_1] I_2 I_1 I_2 \right)
\end{align*}
Similar identities hold for the inverse Murphy elements, and differ only in interchanging $q \leftrightarrow q^{-1}$.
\end{lemma}
\begin{proof}
This follows by a simple, but tedious, application of the relations given in \eqref{eqn:eg_relations}.
\end{proof}
Lemma \ref{lemma:NOddRecursion} and Lemma \ref{lemma:IJ0Iodd} give the corresponding results for odd $N$. As they are proved in a very similar fashion we shall simply state the results.
\begin{lemma}
\label{lemma:NOddRecursion}
For $N$ odd we have:
\begin{align*}
&I_1 J^{(\rm C)}_{2i+1} I_1 = q^{2} I_1 J^{(\rm C)}_{2i} I_1,
&&I_1 J^{(\rm C)}_{2i+2} I_1 = q^{-2} I_1 J^{(\rm C)}_{2i} I_1. 
\end{align*}
The first equality holds for $0 \le i \le (N-3)/2$ and the second
for $0 \le i \le (N-5)/2$. We also have:
\begin{align*}
I_2 J^{(\rm C)}_{2i+2} I_2 = q^2 I_2 J^{(\rm C)}_{2i+1} I_2,\qquad I_2 J^{(\rm C)}_{2i+3} I_2 =
q^{-2} I_2 J^{(\rm C)}_{2i+1} I_2. 
\end{align*}
The first equality holds for $0 \le i \le (N-3)/2$ and the second for $0 \le i \le (N-5)/2$.
\end{lemma}
\begin{lemma}
\label{lemma:IJ0Iodd}
For $N$ odd we have the following identities:
\begin{align*}
I_1 J^{(\rm C)}_0 I_1 &= q^{-2} [2]^{(N-3)/2} \left(
\frac{[\omega_2][2\left(1+\omega_1-\omega_2\right)]}{[\omega_2+1][1+\omega_1-\omega_2]} I_1 \right. \\
& \left. \phantom{\frac{[\omega_2][2\left(1+\omega_1-\omega_2\right)]}{[\omega_2+1][1+\omega_1-\omega_2]} I_1} \quad - (q-q^{-1})^2[1+\omega_1][1-\omega_2]
I_1 I_2 I_1 \right)\\
I_1 J^{(\rm C)}_{N-1} I_1 &= q^{-\omega_2-N+1} [2]^{(N-1)/2}
\left(q^{\omega_1} \frac{[\omega_2]}{[\omega_2+1]} I_1- (q-q^{-1})[1+\omega_1] I_1 I_2 I_1 \right)\\
I_2 J^{(\rm C)}_0 I_2 &= q^{-\omega_1}[2]^{(N-1)/2} \left( q^{\omega_2} 
\frac{[\omega_1]}{[\omega_1+1]} I_2 -(q-q^{-1})  [1+\omega_2]
I_2 I_1 I_2\right)\\
I_2 J^{(\rm C)}_1 I_2 &= q^{-3} [2]^{(N-3)/2} \left( 
\frac{[\omega_1][2 \left(\omega_1-\omega_2 \right)]}{[\omega_1+1][\omega_1-\omega_2]} I_2 + [\omega_1][\omega_2+1] I_2 I_1 I_2  \right)
\end{align*}
Similar identities hold for the inverse Murphy elements,
and differ only in interchanging $q \leftrightarrow q^{-1}$.
\end{lemma}
\subsection{Proof of Theorem~\ref{thm:IJI=bI}}
\begin{proof}[Proof of Theorem \ref{thm:IJI=bI}]
The results of the previous subsection allow us to
evaluate all terms of the form $I_1 J_i^{\pm 1} I_1$ and $I_2 J_i^{\pm 1} I_2$.
\begin{itemize}
\item{$N$ even}

For $N$ even the idempotents given in \eqref{eqn:IdempotentsEven} satisfy:
\begin{align*}
&I_1^2 = [2]^{N/2} I_1,
&&I_2^2 = [2]^{(N-2)/2} \frac{[\omega_1][\omega_2]}{[\omega_1+1][\omega_2+1]} I_2.
\end{align*}
We substitute expressions for $I_1 J_i^{\pm 1} I_1$ into \eqref{eqn:CentreExpanded},  with $I=I_1$, to obtain:
\[
\frac{[2\theta]}{[\theta]} I_1^2 =  [2]^{N/2}
     \left(\frac{[2 \left(\omega_1+\omega_2+1 \right)]}{[\omega_1+\omega_2+1]} I_1 - 
     (q-q^{-1})^2[1+\omega_1][1+\omega_2]I_1I_2I_1\right).
\]
Rearranging this we complete the proof of $I_1I_2I_1 = bI_1$ with $b$
given by \eqref{eqn:beven}. The other case, $I_2I_1I_2=bI_2$, follows similarly.
\item{$N$ odd}

For $N$ odd the idempotents given in \eqref{eqn:IdempotentsOdd} satisfy:
\begin{align*}
&I_1^2 = [2]^{(N-1)/2} \frac{[\omega_2]}{[\omega_2+1]} I_1,
&& I_2^2 = [2]^{(N-1)/2} \frac{[\omega_1]}{[\omega_1+1]} I_2.
\end{align*}
We substitute expressions for $I_1 J_i^{\pm 1} I_1$ into \eqref{eqn:CentreExpanded},  with $I=I_1$, to obtain:
\[
\frac{[2\theta]}{[\theta]} I_1^2 = [2]^{(N-1)/2} 
\left(\frac{[\omega_2][2\left(\omega_1-\omega_2\right)]}{[\omega_2+1][\omega_1-\omega_2]} I_1 + (q-q^{-1})^2[1+\omega_1][\omega_2]I_1I_2I_1 \right).
\]
Rearranging this we complete the proof of $I_1I_2I_1 = bI_1$ with $b$
given by \eqref{eqn:bodd}. The other case, $I_2I_1I_2=bI_2$, follows in a 
similar fashion.
\end{itemize}
\end{proof}
\section{Proof of Proposition~\ref{prop:BulkActionZero}}
\label{ap:BulkActionZero}
This Proposition is proved in two steps. First in Lemma 
\ref{lemma:SimpleBulkActionZero1} and Lemma \ref{lemma:SimpleBulkActionZero2} we prove that the action of the bulk generators vanishes on some simple slopes. When acting with a bulk generator on a more general slope we use the Yang-Baxter and right reflection equations to reduce the problem to simpler slopes.
\subsection{Simple slopes}
\begin{lemma}
\label{lemma:SimpleBulkActionZero1}
The following identities hold in the 2BTL algebra:
\begin{align}
\label{eqn:BulkFundIdentities1a}
e_{2n} R_{2n \pm 1}(\omega_1) E_N&=0 \quad \quad 1 \le 2n, 2n \pm 1 \le N-1 \\
e_{2n \pm 1} R_{2n}(-\omega_1-1) E_N &= 0 \nonumber
\end{align}
\end{lemma}
\begin{proof}
It is sufficient to prove the first line of \eqref{eqn:BulkFundIdentities1a} as the second follows by applying $e_{2n \pm 1}$ and using the identity:
\[
e_{2n \pm 1} e_{2n} R_{2n \pm 1}(u)  = - \f{[u+1]}{[u]} e_{2n \pm 1} R_{2n} (-u-1) 
\]
We now prove the first line of \eqref{eqn:BulkFundIdentities1a}. Here we only give the $+$ case as the other is similar:
\begin{align*}
e_{2n} R_{2n+1}(\omega_1) E_{2n+2} &= \f{[\omega_1+1]}{[\omega_1]} E_{2n-1} e_{2n} R_{2n+1}(\omega_1) e_{2n-1} e_{2n} e_{2n+1} E_{2n+1} \\
&= \f{[\omega_1+1]}{[\omega_1]} e_{2n} e_{2n+1} E_{2n-1} \left( e_{2n-1}-\f{[\omega_1+1]}{[\omega_1]} \right) E_{2n+1}  \\
&= 0
\end{align*}
\end{proof}
Recalling the (un-normalized) idempotents $I_1$ and $I_2$ which are
defined in \eqref{eqn:IdempotentsEven} and \eqref{eqn:IdempotentsOdd}, we have:
\begin{lemma}
\label{lemma:SimpleResult}
The following identities hold in the 2BTL algebra:
\begin{itemize}
\item For $N$ even:
\begin{align*}
&I_2 I_1 E_N = s_1^{-N/2} I_2 E_N ,
&& I_1 I_2 E_N = s_1^{N/2} I_1 e_N E_N
\end{align*}
\item For $N$ odd:
\begin{align*}
&I_1 I_2 E_N = s_1^{(N+1)/2} I_1 E_N,
&& I_2 I_1 E_N = s_1^{-(N-1)/2} I_2 e_N E_N
\end{align*}
\end{itemize}
where $s_1=\f{[\omega_1]}{[\omega_1+1]}$.
\end{lemma}
\begin{proof}
The proofs are similar in both cases and we shall only give $N$
even. We define:
\begin{align*}
&I_{1,i} = e_{2i+1} e_{2i+3} \cdots e_{N-1},
&&I_{2,i} = e_{2i} e_{2i+2} \cdots e_{N}
\end{align*}
Using Lemma \ref{lemma:e0action} and Lemma \ref{lemma:SimpleBulkActionZero1} we have:
\begin{align*}
&I_2 I_{1,i} E_N = s_1^{-1} I_2 I_{1,i+1} E_N 
&&I_1 I_{2,i} E_N = s_1 I_1 I_{2,i+1} E_N
\end{align*}
and the Proposition follows.
\end{proof}
\begin{lemma}
\label{lemma:SimpleBulkActionZero2}
For the 2BTL algebra in the double quotient we have the following identities:
\begin{itemize}
\item For $N$ even:
\begin{align*}
&e_{N-1} K_N(-\omega_1-1) E_N =0, 
&&e_{N-1} K_N(\omega_1-1) R_{N-1}(\omega_1) E_N = 0
\end{align*}
\item For odd $N$:
\begin{align*} 
&e_{N-1} K_N(\omega_1) E_N =0, 
&&e_{N-1} K_N(-\omega_1-2) R_{N-1}(-\omega_1-1) E_N = 0
\end{align*}
\end{itemize}
\end{lemma}
\begin{proof}
It is sufficient to prove only the first relation in each case due to the identity:
\[
e_{N-1} K_N(u-1) R_{N-1}(u) = - \f{[u+1]}{[u]} e_{N-1} K_N(-u-1). 
\]
Consider $N$ even. As in the proof of Lemma~\ref{lemma:SimpleResult} we define:
\[
I_{1,i} = e_{2i+1} e_{2i+3} \cdots e_{N-1}
\]
The identity follows as a special case, namely $i=N/2-1$, of the following identity:
\begin{align}
\label{eqn:MoreGeneralIdentity}
I_{1,i} K_N(-\omega_1-1) E_N=0 \qquad 0\leq i \leq N/2-1,
\end{align}
We shall prove this inductively. From Lemma \ref{lemma:SimpleResult} we have:
\[
I_1 I_2 I_1 E_N = I_1 e_N E_N
\]
Now using the double quotient and the identity $b=k(-\omega_1-1)$ we prove the $i=0$ case of \eqref{eqn:MoreGeneralIdentity}. 

We now assume that \eqref{eqn:MoreGeneralIdentity} holds for some $i=n$ with $0 \le n \le N/2-2$. We have:
\begin{align*}
I_{1,n} K_N(-\omega_1-1) E_N &= e_{2n+1} I_{1,n+1} K_N(-\omega_1-1) E_{2n+1} E_N \\
&= e_{2n+1} E_{2n+1} I_{1,n+1} K_N(-\omega_1-1) E_N
\end{align*}
where we have used the definition of $I_{1,n}$ and Proposition \ref{lemma:1BTLStartingVectors}. The case $i=n+1$ of \eqref{eqn:MoreGeneralIdentity} now follows from pre-multiplying by $E_{2n+1}$, using Definition \ref{defn:1BTLStartingVectors}, and $\left[ e_N, E_{2n+2} \right]=0$.

For $N$ odd the proof is similar and we instead use the identity $k(\omega_1) = b [\omega_1+1]/[\omega_1]$.  
\end{proof}

The results of Lemma \ref{lemma:SimpleBulkActionZero1} can be
interpreted in a pictorial way. They correspond to vanishing of the
bulk generators on some simple slopes. For example $R_{2n-1}(\omega_1)
E_N$ corresponds to a single tile being added to the fundamental path
at position $2n-1$. From the form of the fundamental path we know that
this tile must be added from above: 
\begin{align*}
{\includegraphics[height=3 cm]{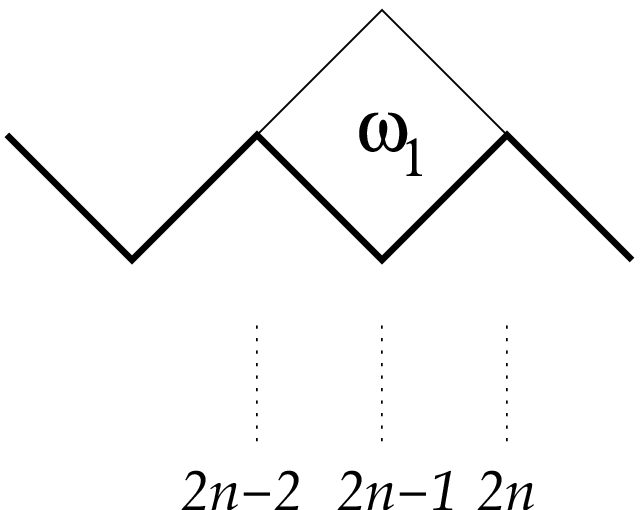}}
\end{align*}
The results of Lemma \ref{lemma:SimpleBulkActionZero1} now imply that the generators $e_{2n}$ and $e_{2n-2}$ vanish on this one-tile slope. We have similar pictures for slopes created by a single tile added from below.

In a similar way we can interpret the results of Lemma  \ref{lemma:SimpleBulkActionZero2}. For $N$ even we have vanishing of the bulk generator $e_{N-1}$ on the two simple slopes:  
\begin{align*}
&\includegraphics[width=5 cm]{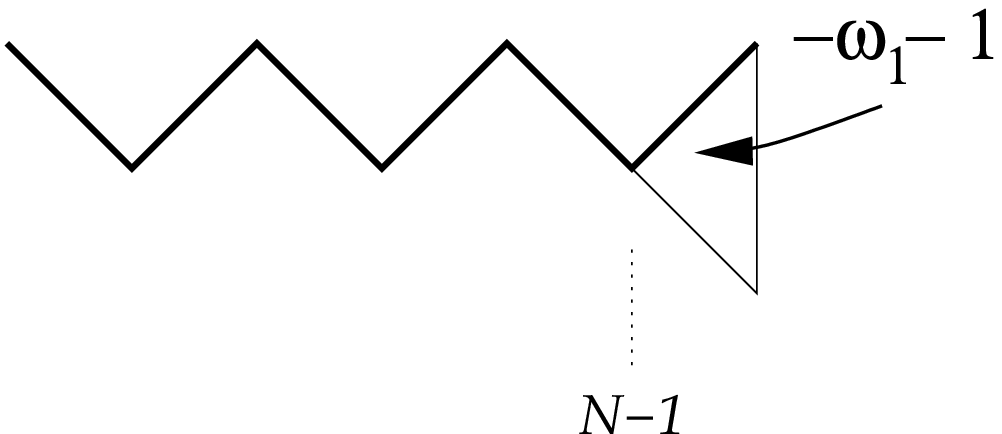} 
&\includegraphics[width=5 cm]{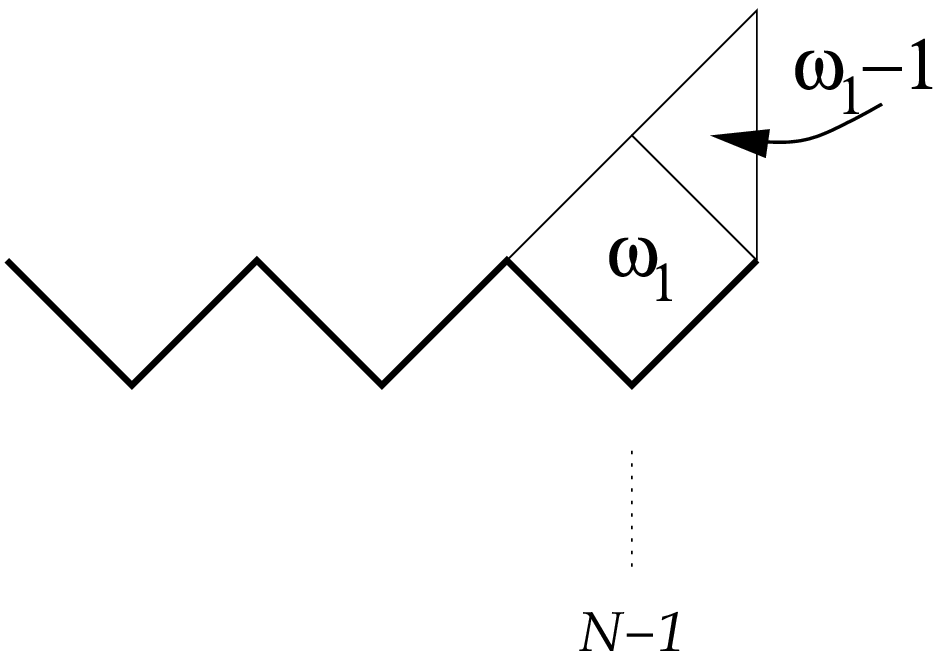}
\end{align*}
There is a similar set of pictures for $N$ odd.
\subsection{Proof of Proposition~\ref{prop:BulkActionZero}}
\begin{proof}[Proof of Proposition \ref{prop:BulkActionZero}]
We use the Yang-Baxter and right reflection equations to move the bulk
operators $e_i=R_{i}(-1)$ through expressions until we reduce
ourselves to the cases covered in Lemma
\ref{lemma:SimpleBulkActionZero1} and Lemma \ref{lemma:SimpleBulkActionZero2}.

We shall use white boxes to denote $R_{i}(u)$ or $K_N(u)$ and grey ones to represent
bulk $e_i=R_i(-1)$ generators. We shall only consider the case $h \ge0$ as the
other case $h<0$ is similar. 
\begin{itemize}
\item{$e_i$ acting on a descending slope with $h\ge0$}

Writing $u_h=\omega_1-h$ the action of $e_i=R_i(-1)$ on a slope of a
path is locally depicted by the left hand side of
Figure~\ref{fig:eiOnSlope}. Using the Yang-Baxter equation:
\[
R_{i}(-1)R_{i-1}(u_h)R_i(u_h+1) = R_{i-1}(u_h+1)R_i(u_h) R_{i-1}(-1) 
\]
we can pull the grey box through and note that the generator
$R_{i-1}(-1)$ again will act on a descending 
slope. We now repeat until eventually we reach a point at which we can
use Lemma  \ref{lemma:SimpleBulkActionZero1} or Lemma \ref{lemma:SimpleBulkActionZero2}.
\begin{figure}[h]
\centerline{\includegraphics[width=8 cm]{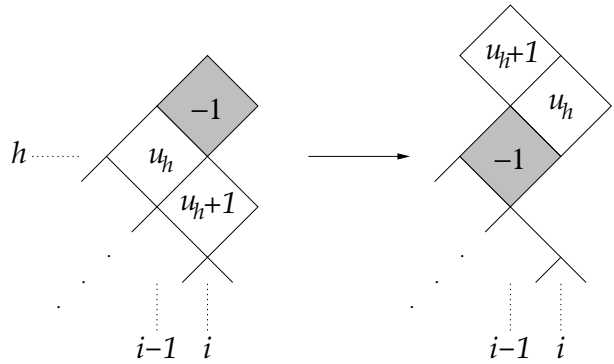}}
\caption{Action of $e_i$ on a slope with $h \ge 0$. Due to the Yang-Baxter
equation we can pull through the grey tile.}
\label{fig:eiOnSlope}
\end{figure}
\item{$e_i$ acting on an ascending slope with $h \ge 0$}
\begin{align*}
\includegraphics[width=4 cm]{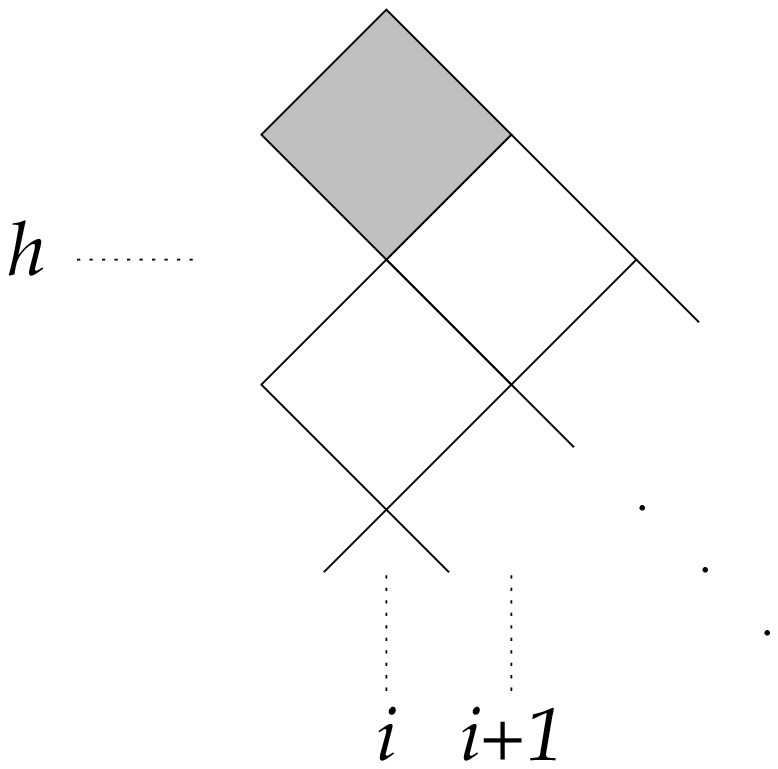}
\end{align*}
Once again we use the Yang-Baxter equation as before to move $e_i$
through. There are now two possibilities: we will reach a point at
which we can use Lemma \ref{lemma:SimpleBulkActionZero1} or
Lemma \ref{lemma:SimpleBulkActionZero2} we will reach the right boundary. In the latter case we have:
\begin{align*}
\includegraphics[width=3 cm]{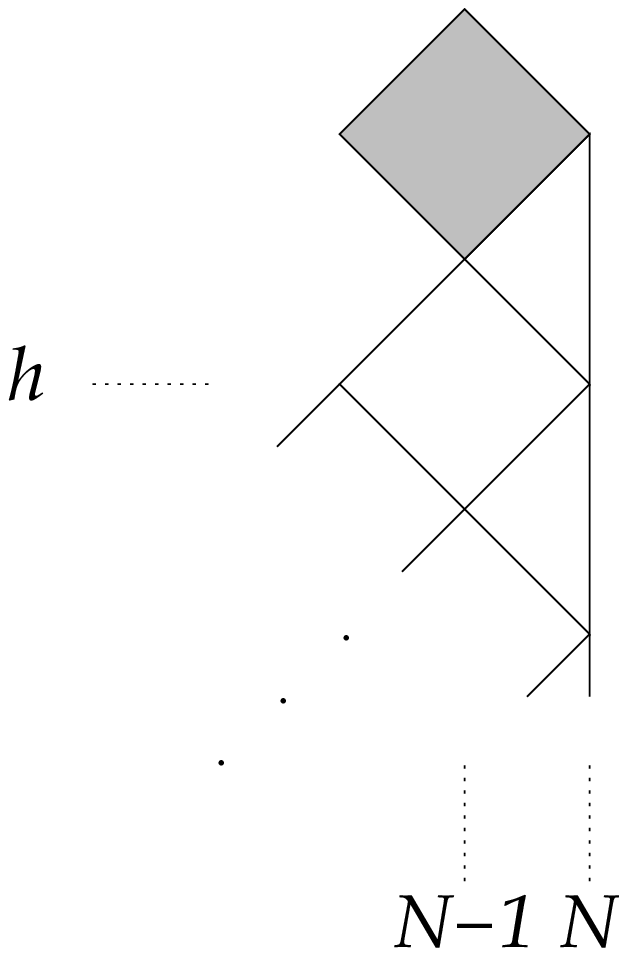}
\end{align*}
We now use the right reflection equation, with $u_h=\omega_1-h$,:
\begin{multline*}
R_{N-1}(-1) K_N(u_h-1) R_{N-1}(u_h) K_N(u_h+1) \\
= K_N(u_h+1)R_{N-1}(u_h)K_N(u_h-1)R_{N-1}(-1)
\end{multline*}
The generator $e_{N-1}$ now acts on descending paths with $h \ge 0$ and we can use the previous result to show this vanishes.
\end{itemize}
\end{proof}
%
%



\providecommand{\href}[2]{#2}\begingroup\raggedright
\endgroup

\end{document}